\documentclass[aop,seceqn,MSNbibl,citesort,noMRlinks,dvips]{arximspdf}

\doi{10.1214/08-AOP394}
\volume{37}
\issue{1}
\pubyear{2009}
\firstpage{1}
\lastpage{47}

\makeatletter

\newtheorem{theorem}{Theorem}[section]
\newproclaim{remark}{Remark}[section]
\newtheorem{lemma}{Lemma}[section]
\newtheorem{proposition}{Proposition}[section]
\newproclaim{step}{Step}

\makeatother

\begin{document}
\begin{frontmatter}

\title{The largest eigenvalues of finite rank deformation of large Wigner
matrices: convergence and nonuniversality of the fluctuations}
\runtitle{Largest eigenvalues of deformed Wigner
matrices}

\begin{aug}
\author[A]{\fnms{Mireille} \snm{Capitaine}\ead[label=e1]{mireille.capitaine@math.univ-toulouse.fr}},
\author[B]{\fnms{Catherine} \snm{Donati-Martin}\ead[label=e2]{catherine.donati@upmc.fr}}\\
and
\author[C]{\fnms{Delphine} \snm{F\'eral}\corref{}\ead[label=e3]{delphine.feral@math.u-bordeaux1.fr}}
\runauthor{M. Capitaine, C. Donati-Martin and D. F\'eral}
\affiliation{Institut de Math\'{e}matiques de Toulouse, Universit\'{e} Paris 6
and~Instituts~de~Math\'{e}matiques de Toulouse et de Bordeaux}
\address[A]{M. Capitaine\\
Institut de Math\'ematiques de Toulouse\\
Equipe de Statistique et Probabilit\'es\\
31062 Toulouse Cedex 09\\
France \\
\printead{e1}}
\address[B]{C. Donati-Martin\\
UPMC Universit\'{e} Paris 6\\
and\\
CNRS, UMR 7599\\
Probabilit\'es et Mod\`eles Al\'eatoires\\
Site Chevaleret\\
16 rue Clisson\\
75013 Paris\\
France\\
\printead{e2}}
\address[C]{D. F\'eral\\
Institut de Math\'{e}matiques de Bordeaux\\
Universit\'{e} Bordeaux 1\\
351 Cours de la Lib\'eratation\\
33405 Talence Cedex\\
France\\
\printead{e3}}
\end{aug}
\pdfauthor{Mireille Capitaine, Catherine  Donati-Martin, Delphine Feral}

\received{\smonth{5} \syear{2007}}
\revised{\smonth{12} \syear{2007}}

%
\begin{abstract}
In this paper, we investigate the asymptotic spectrum of complex or
real Deformed Wigner matrices $(M_N)_N$ defined by $M_N= W_N /
\sqrt{N} + A_N$ where~$W_N$ is an $N \times N$ Hermitian (resp.,
symmetric) Wigner matrix whose entries have a symmetric law
satisfying a Poincar\'e inequality. The matrix $A_N$ is Hermitian
(resp., symmetric) and deterministic with all but finitely many
eigenvalues equal to zero. We first show that, as soon as the first
largest or last smallest eigenvalues of $A_N$ are sufficiently far
from zero, the corresponding eigenvalues of $M_N$ almost surely exit
the limiting semicircle compact support as the size $N$ becomes
large. The corresponding limits are universal in the sense that they
only involve the variance of the entries of $W_N$. On the other
hand, when $A_N$ is diagonal with a sole simple nonnull eigenvalue
large enough, we prove that the fluctuations of the largest
eigenvalue are not universal and vary with the particular
distribution of the entries of $W_N$.
\end{abstract}

%
\pdfkeywords{Deformed Wigner matrices,
asymptotic spectrum,
Stieltjes transform,
largest eigenvalues,
fluctuations,
nonuniversality}

\end{frontmatter}

\section{Introduction}\label{intro}
This paper lies in the lineage of recent works studying the
influence of some perturbations on the asymptotic spectrum of
classical random matrix models. Such questions come from statistics
(cf. \cite{John}) and appeared in the framework of empirical
covariance matrices, also called nonwhite Wishart matrices or
spiked population models, considered by Baik, Ben Arous and
P\'ech\'e \cite{BBP} and by Baik and Silverstein \cite{BS3}.
The work \cite{BBP} deals with random sample covariance matrices
$(S_N)_N$ defined by
%
\begin{eqnarray}\label{spike}
S_N=\frac{1}{N} Y_N^* Y_N,
\end{eqnarray}
where $Y_N$ is a $p \times N$ complex matrix whose sample column
vectors are i.i.d., centered, Gaussian and of covariance matrix a
deterministic Hermitian matrix ${\Sigma}_p$ having all but finitely
many eigenvalues equal to 1. Besides, the size of the samples $N$
and the size of the population $p=p_N$ are assumed of the same order
(as $N \to\infty$). The authors of \cite{BBP} first noticed that,
as in the classical case (known as the Wishart model) where
$\Sigma_p=I_p$ is the identity matrix, the global limiting behavior
of the spectrum of $S_N$ is not affected by the matrix $\Sigma_p$.
Thus, the limiting spectral measure is the well-known
Marchenko--Pastur law. On the other hand, they pointed out a phase
transition phenomenon for the fluctuations of the largest eigenvalue
according to the value of the largest eigenvalue(s) of $\Sigma_p$.
The approach of \cite{BBP} does not extend to the real Gaussian
setting and the whole analogue of their result is still an open
question. Nevertheless, Paul was able to establish in \cite{Pa}
the Gaussian fluctuations of the largest eigenvalue of the real
Gaussian matrix $S_N$ when the largest eigenvalue of $\Sigma_p$ is
simple and sufficiently larger than 1. More recently, Baik and
Silverstein investigated in \cite{BS3} the almost sure limiting
behavior of the extremal eigenvalues of complex or real
nonnecessarily Gaussian matrices. Under assumptions on the first four
moments of the entries of $Y_N$, they showed in particular that when
exactly $k$ eigenvalues of ${\Sigma}_p$ are far from 1, the $k$
first eigenvalues of $S_N$ are almost surely outside the limiting
Marchenko--Pastur support. Fluctuations of the eigenvalues that jump
are universal and have been recently found by Bai and Yao
in \cite{BY2} (we refer the reader to \cite{BY2} for the precise
restrictions made on the definition of the covariance matrix $\Sigma
_p$). Note that the problem of the fluctuations in the
very general setting of \cite{BS3} is still open.

Our purpose here is to investigate the asymptotic behavior of the
first extremal eigenvalues of some complex or real Deformed Wigner
matrices. These models can be seen as the additive analogue of the
spiked population models and are defined by a sequence $(M_N)_N$
given by
%
\begin{eqnarray} {\label{Dw}}
M_N=\frac{1}{\sqrt{N}} W_N +A_N:= X_N + A_N,
\end{eqnarray}
where $W_N$ is a Wigner matrix such that the common distribution of
its entries satisfies some technical conditions [given in (i)
below] and $A_N$ is a deterministic matrix of finite rank. We
establish the analogue of the main result of \cite{BS3}, namely that,
once~$A_N$ has exactly $k$ (fixed) eigenvalues far enough from zero,
the $k$ first eigenvalues of $M_N$ jump almost surely outside the
limiting semicircle support. This result is universal (as the one of
\cite{BS3}) since the corresponding limits only involve the variance
of the entries of $W_N$. On the other hand, at the level of the
fluctuations, we exhibit a striking phenomenon in the particular
case where $A_N$ is diagonal with a sole simple nonnull eigenvalue
large enough. Indeed, we find that in this case, the fluctuations of
the largest eigenvalue of $M_N$ are not universal and strongly
depend on the particular law of the entries of $W_N$. More
precisely, we prove that the limiting distribution of the (properly
rescaled) largest eigenvalue of~$M_N$ is the convolution of the
distribution of the entries of $W_N$ with a Gaussian law. In
particular, if the entries of $W_N$ are not Gaussian,
the fluctuations of the largest eigenvalue of $M_N$ are not Gaussian.

In the following section, we first give the precise definition of
the Deformed Wigner matrices (\ref{Dw}) considered in this paper and
we recall the known results on their asymptotic spectrum. Then, we
present our results and sketch the proofs. We also outline the
organization of the paper.

\section{Model and results}\label{results}
Throughout this paper, we consider complex or real Deformed Wigner
matrices $(M_N)_N$ of the form (\ref{Dw}) where the matrices $W_N$
and~$A_N$ are defined as follows:
\begin{itemize}[(ii)]
\item[(i)] $W_N$ is an $N \times N$ Wigner Hermitian (resp.,
symmetric) matrix such that the~$N^2$ random variables
$(W_N)_{ii}$, $\sqrt{2} \Re e((W_N)_{ij})_{i<j}$, $\sqrt{2} \Im m
((W_N)_{ij})_{i<j}$ [resp., the $N(N+1)/2$ random variables
$\frac{1}{\sqrt{2}} (W_N)_{ii}$, $(W_N)_{ij}$, $i<j$] are
independent identically distributed with a symmetric distribution
$\mu$ of variance $\sigma^2$ and satisfying a Poincar\'e inequality
(see Section \ref{base}).
\item[(ii)] $A_N$ is a deterministic Hermitian (resp., symmetric)
matrix of fixed finite
rank~$r$ and built from a family of $J$ fixed real numbers
$\theta_1>\cdots> \theta_J$ independent of~$N$ with some $j_0$ such
that $\theta_{j_0}=0$. We assume that the nonnull eigenvalues~$\theta_j$
of $A_N$ are of fixed multiplicity $k_{j}$ (with $\sum_{j
\not= j_0} k_j=r$), that is, $A_N$ is similar to the diagonal matrix
\end{itemize}
%
%
\begin{equation}\label{defdiagAN}\qquad\quad
D_N = \operatorname{diag}\biggl(\underbrace{\theta_1,\ldots
,\theta_1}_{k_{1}}, \ldots,\underbrace{\theta_{j_0-1}, \ldots
}_{k_{{j_0-1}}} \,  ,\underbrace{0, \ldotp\ldots, 0}_{N-r},
\underbrace{\theta_{j_0+1},\ldots
}_{k_{{j_0+1}}}  \, , \ldots, \underbrace{\theta_{J},\ldots
}_{k_J}\biggr).
\end{equation}

Before going into the details of the results, we want to point out
that the condition made on $\mu$ (namely that $\mu$ satisfies a
Poincar\'e inequality) is just a technical condition: we conjecture
that our results still hold under weaker assumptions (see Remark
\ref{Rmk1} below). Nevertheless, a lot of measures satisfy a
Poincar\'e inequality (we refer the reader to \cite{B} for a
characterization of such measures on $\mathbb R$; see also~\cite
{Tou}). For
instance,
consider $\mu(dx) = \exp(-|x|^\alpha) \,dx$ with $\alpha\geq1$.

Furthermore, note that this condition implies that $\mu$ has moments
of any order (cf.
Corollary 3.2 and Proposition 1.10 in \cite{L}).

Let us now introduce some notations. When the entries of $W_N$ are
further assumed to be Gaussian, that is, in the complex (resp., real)
setting when $W_N$ is of the so-called GUE (resp., GOE), we will
write $W_N^G$ instead of $W_N$. Then $X_N^G := W_N^G/ \sqrt N$ will
be said to be of the GU(O)E($N,
\frac{\sigma^2}{N}$) and we will let $M_N^G=X_N^G+A_N$ be the
corresponding Deformed GU(O)E model.

In the following, given an arbitrary Hermitian matrix $B$ of order
$N$, we will denote by $\lambda_1(B) \geq\cdots\geq\lambda_N(B)$
its $N$ ordered eigenvalues and by $\mu_B=\frac{1}{N} \sum_{i=1}^N
\delta_{\lambda_i(B)}$ its empirical measure. $\operatorname{Spect}(B)$ will
denote the spectrum of $B$. For notational convenience, we
will also set $\lambda_0(B)= + \infty$ and $\lambda_{N+1}(B)= -
\infty
$.

The Deformed Wigner model is built in such a way that the Wigner
theorem is still satisfied. Thus, as in the classical Wigner model
($A_N \equiv0$), the spectral measure $(\mu_{M_N})$ converges a.s.
toward the semicircle law $\mu_{sc}$ whose density is given by
%
\begin{equation}\label{scl}
\frac{d \mu_{sc}}{dx}(x)= \frac{1}{2 \pi\sigma^2} \sqrt{4\sigma^2
- x^2}   1_{[-2\sigma,2 \sigma]}(x).
\end{equation}
This result follows from Lemma 2.2 of \cite{Bai}. Note that it only
relies on the two first moment assumptions on the entries of $W_N$
and the
fact that the $A_N$'s are of finite rank.

On the other hand, the asymptotic behavior of the extremal
eigenvalues may be affected by the perturbation $A_N$. Recently,
P\'ech\'e studied in \cite{Pe} the Deformed GUE under a finite rank
perturbation $A_N$ defined by (ii). Following the method of~\cite{BBP},
she highlighted the effects of the nonnull eigenvalues
of $A_N$ at the level of the fluctuations of the largest eigenvalue
of $M_N^G$. To explain this in more detail, let us recall that when
$A_N \equiv0$, it was established in \cite{TW} that as $N
\rightarrow\infty$,
%
\begin{equation}\label{fluctuGUE}
\sigma^{-1} N^{2/3} \bigl( \lambda_1(X_N^G) - 2 \sigma\bigr)
\stackrel{\mathcal L}{\longrightarrow} F_2,
\end{equation}
where $F_2$ is the well-known GUE Tracy--Widom distribution (see
\cite{TW} for the precise definition). Dealing with the Deformed GUE
$M_N^G$, it appears that this result is modified as soon as the
first largest eigenvalue(s) of $A_N$ is (are) quite far from zero.
In the particular case of a rank-1 perturbation $A_N$ having a
fixed nonnull eigenvalue $\theta>0$, \cite{Pe} proved that the
fluctuations of the largest eigenvalue of $M_N^G$ are still given by
(\ref{fluctuGUE}) when $\theta$ is small enough and precisely when
$\theta<\sigma$. The limiting law is changed when $\theta=\sigma$.
As soon as $\theta>\sigma$, \cite{Pe} established that the largest
eigenvalue $\lambda_1(M_N^G)$ fluctuates around
%
\begin{equation}\label{rho}
\rho_{\theta}=\theta+\frac{\sigma^2}{\theta}
\end{equation}
(which is $>2\sigma$ since $\theta> \sigma$) as
%
\begin{equation}\label{PecheDefGUE}
\sqrt{N} \bigl(
\lambda_{1}(M_N^G)-\rho_{\theta}\bigr) \stackrel{\mathcal
L}{\longrightarrow}
\mathcal N (0,\sigma_{\theta} ^2),
\end{equation}
where
%
\begin{equation}\label{sigma}
\sigma_{\theta}=(\sigma/\theta) \sqrt{{\theta^2} - {\sigma^2}}.
\end{equation}
Similar results are conjectured for the Deformed GOE but
P\'ech\'e emphasized that her approach fails in the real framework.
Indeed, it is based on the explicit Fredholm determinantal
representation for the distribution of the largest eigenvalue(s)
that is specific to the complex setting. Nevertheless, Ma\"ida
\cite{Ma} obtained a large deviation principle for the largest
eigenvalue of the Deformed GOE $M_N^G$ under a rank-1 deformation
$A_N$; from this result she could deduce the almost sure limit with
respect to the nonnull eigenvalue of $A_N$. Thus, under a rank-1
perturbation~$A_N$ such that $D_N=\operatorname{diag}(\theta, 0, \ldots, 0)$
where $\theta>0$, \cite{Ma} showed that
%
\begin{equation}\label{asGUE1}
\lambda_{1}(M_N^G)
\stackrel{\mathrm{a.s.}}{\longrightarrow} \rho_{\theta}\qquad \mbox{if }
\theta
> \sigma
\end{equation}
and
%
\begin{equation}\label{asGUE2}
\lambda_{1}(M_N^G)
\stackrel{\mathrm{a.s.}}{\longrightarrow} 2 \sigma  \qquad\mbox{if } \theta
\leq
\sigma.
\end{equation}
Note that the approach of \cite{Ma} extends with minor modifications
to the Deformed GUE. Following the investigations of \cite{BS3} in
the context of general spiked population models, one can conjecture
that such a phenomenon holds in a more general and nonnecessarily
Gaussian setting. The first result of our paper, namely the
following Theorem \ref{ThmASCV}, is related to this question. Before
being more explicit, let us recall that when $A_N \equiv0$, the
whole spectrum of the rescaled complex or real Wigner matrix
$X_N=W_N/{\sqrt{N}}$ belongs almost surely to the semicircle support
$[-2\sigma,2\sigma]$ as $N$ goes to infinity and that (cf.
\cite{BYi} or Theorem 2.12 in \cite{Bai})
%
\begin{equation} {\label{extremal}}
\lambda_1(X_N) \stackrel{\mathrm{a.s.}}{\longrightarrow} 2 \sigma
\quad\mbox{and} \quad \lambda_N(X_N) \stackrel{\mathrm{a.s.}}{\longrightarrow}-2
\sigma.
\end{equation}
Note that this last result holds true in a more general setting than
the one considered here (see \cite{BYi} for details) and in
particular only requires the finiteness of the fourth moment of the
law $\mu$. Moreover, one can readily extend the previous limits to
the first extremal eigenvalues of $X_N$, that is,
%
\begin{equation}\label{extremalbis}\qquad\quad
\mbox{for any fixed $k \geq
1$,} \qquad     \lambda_k(X_N) \stackrel{\mathrm{a.s.}}{\longrightarrow} 2
\sigma \quad\mbox{and}\quad  \lambda_{N-k}(X_N)
\stackrel{\mathrm{a.s.}}{\longrightarrow}-2 \sigma.
\end{equation}
Here, we
prove that, under the assumptions (i)--(ii),
(\ref{extremalbis}) fails when some of the $\theta_j$'s are
sufficiently far from zero: as soon as some of the first largest
(resp., last smallest) nonnull eigenvalues $\theta_j$ of $A_N$ are
taken strictly larger than $\sigma$ (resp., strictly smaller than
$-\sigma$), the same part of the spectrum of $M_N$ almost surely
exits the semicircle support $[-2\sigma,2 \sigma]$ as $N \to\infty$
and the new limits are the $\rho_{\theta_j}$'s defined by
%
\begin{equation}{\label{defrhotheta}}
\rho_{\theta_j}=\theta_j + \frac{\sigma^2}{\theta_j}.
\end{equation}
Observe that $\rho_{\theta_j}$ is $>2 \sigma$ (resp., $<-2 \sigma$)
when $\theta_j> \sigma$ (resp., $<-\sigma$) (and $\rho
_{\theta_j} =\pm2 \sigma$ if $\theta_j= \pm\sigma$).

Here is the precise formulation of our result. For definiteness, we
set $k_1+ \cdots+k_{j-1}:=0$ if $j=1$.
\begin{theorem}{\label{ThmASCV}}
Let $J_{+ \sigma}$ (resp., $J_{- \sigma}$) be the number of j's such
that $\theta_j > \sigma$ (resp., $\theta_j < -\sigma$).
\begin{longlist}[(a)]
\item[(a)] $ \forall1 \leq j \leq J_{+ \sigma},   \forall1
\leq i \leq k_j,  \lambda_{k_1+ \cdots+k_{j-1} +i}(M_N)
\longrightarrow\rho_{\theta_j}  \mbox{ a.s.}$
\item[(b)] $ \lambda_{k_1+ \cdots+k_{J_{+ \sigma}} +1}(M_N)
\longrightarrow2\sigma \mbox{ a.s.}$
\item[(c)] $ \lambda_{k_1+ \cdots+k_{J-J_{- \sigma}}}(M_N)
\longrightarrow-2 \sigma \mbox{ a.s.}$
\item[(d)] $ \forall j \geq J-J_{- \sigma}+1,   \forall1
\leq i \leq k_j,
\lambda_{k_1+ \cdots+k_{j-1} +i}(M_N) \longrightarrow\rho_{\theta_j}
 \mbox{ a.s}.$
\end{longlist}
\end{theorem}

\begin{remark}{\label{Rmk1}} Following \cite{BS3}, one can expect that
this theorem holds true in a more general setting than the one considered
here, namely one that would only require four
first moment conditions on the law $\mu$ of the Wigner entries. As we will
explain in the following, the assumption that $\mu$ satisfies a
Poincar\'e inequality is actually fundamental in our
reasoning since we will need several variance estimates.
\end{remark}

This theorem will be proved in Section \ref{section-as}. The second
part of this work is devoted to the study of the particular rank-1
diagonal deformation $A_N=\operatorname{diag}(\theta, 0, \ldots,0)$ such
that $\theta> \sigma$. We investigate the fluctuations of the
largest eigenvalue of any real or complex Deformed model $M_N$
satisfying (i) around its limit $\rho_\theta$. We obtain the
following result.

\begin{theorem}{\label{ThmFluctuations}}
Let $A_N=\operatorname{diag}(\theta,0,\ldots, 0)$ with $\theta>
\sigma$. Define
%
\begin{equation}\label{defvtheta}
v_{\theta}= \frac{t}{4}\biggl( \frac{m_4 - 3 \sigma^4}{\theta^2}
\biggr) + \frac{t}{2}\frac{\sigma^4}{\theta^2-\sigma^2},
\end{equation}
where $t=4$ (resp., $t=2$) when $W_N$ is real (resp., complex) and
$m_4:= \int x^4\, d\mu(x)$. Then
%
\begin{equation}\label{FluctDef}
\sqrt{N} \biggl(1- \frac{\sigma^2}{\theta^2} \biggr)^{-1} \bigl(\lambda_1(M_N)
-\rho_{\theta} \bigr) \stackrel{\mathcal
L}{\longrightarrow} \mu\ast\mathcal{N}(0, v_{\theta}).
\end{equation}
\end{theorem}

Note that when $m_4 = 3 \sigma^4$ as in the Gaussian case,
the variance of the limiting distribution of $\sqrt{N} (\lambda_1(M_N)
-\rho_{\theta} )$ is equal
to
$\sigma_{\theta} ^2$ (resp., $2\sigma_{\theta} ^2$) in the complex
(resp.,
real) setting [with $\sigma_{\theta}$ given by (\ref{sigma})].

\begin{remark}
Since $\mu$ is symmetric, it readily follows from Theorem
\ref{ThmFluctuations} that when $A_N=\operatorname{diag}(\theta,0,\ldots,
0)$ and $\theta<- \sigma$, the smallest eigenvalue of $M_N$
fluctuates as $ \sqrt{N} (1- \sigma^2/\theta^2 )^{-1} (\lambda_N(M_N)
-\rho_{\theta} ) \stackrel{\mathcal
L}{\longrightarrow} \mu\ast\mathcal{N}(0, v_{\theta}).$
\end{remark}

In particular, one derives the analogue of
(\ref{PecheDefGUE}) for the Deformed GOE:

\begin{theorem}{\label{ThmFluctDefGOE}}
Let $A_N$ be an arbitrary deterministic symmetric matrix of rank 1
having a nonnull eigenvalue $\theta$ such that $\theta> \sigma$.
Then the largest eigenvalue of the Deformed GOE fluctuates as
%
\begin{equation}\label{FluctDefGOE}
\sqrt{N} \bigl(\lambda_1(M_N^G)-\rho_{\theta} \bigr) \stackrel
{\mathcal
L}{\longrightarrow} \mathcal N ( 0, 2 \sigma_{\theta} ^2).
\end{equation}
\end{theorem}

Obviously, thanks to the orthogonal invariance of the GOE,
this result is
a direct consequence of Theorem \ref{ThmFluctuations}.

It is worth noticing that, according to the Cram\'{e}r--L\'evy theorem
(cf. \cite{Fel}, Theorem 1, page 525), the limiting distribution
(\ref{FluctDef}) is not Gaussian if $\mu$ is not Gaussian. Thus,
(\ref{FluctDef}) depends on the particular law $\mu$ of the entries
of the Wigner matrix $W_N$ which implies the nonuniversality of the
fluctuations of the largest eigenvalue of rank-1 diagonal
deformation of symmetric or Hermitian Wigner matrices (as
conjectured in Remark 1.7 of \cite{FePe}).

The latter also shows that in the non-Gaussian setting, the
fluctuations of the largest eigenvalue depend, not only on the
spectrum of the deformation $A_N$, but also on the particular
definition of the matrix $A_N$. Indeed, in collaboration with S.
P\'ech\'e, the third author of the present article has recently
stated in \cite{FePe} the universality of the fluctuations of some
Deformed Wigner models under a full deformation $A_N$ defined by
$(A_N)_{ij} = {\theta}/{N}$ for all $1\leq i,j \leq N$ (see also
\cite{FK}). Before giving some details on this work, we have to
specify that \cite{FePe} considered Deformed models such that the
entries of the Wigner matrix $W_N$ have sub-Gaussian moments.
Nevertheless, thanks to the analysis made in \cite{Ru}, one can
observe that the assumptions of \cite{FePe} can be reduced and that
it is, for example, sufficient to assume that the $W_{i,j}$'s have
moments of any order. Thus, the conclusions of \cite{FePe} apply to
the setting considered in our paper. The main result of \cite{FePe}
establishes the universality of the fluctuations of the largest
eigenvalue of the complex Deformed model $M_N$ associated to a full
deformation $A_N$ and for any value of the parameter~$\theta$. In
particular, when $\theta>\sigma$, it is proved therein the
universality of the Gaussian fluctuations (\ref{PecheDefGUE}). The
approach of \cite{FePe} is mainly based on a combinatorial method
inspired by the work \cite{So} (which handles the non-Deformed
Wigner model) and some results of \cite{Pe} on the Deformed GUE. The
combinatorial arguments of \cite{FePe} also work (with minor
modifications) in the real framework and yield the universality of
the fluctuations if $\theta< \sigma$. In the case where $\theta>
\sigma$ which is of particular interest here, the analysis made in
\cite{FePe} reduces the universality problem in the real setting to
the knowledge of the particular Deformed GOE model (this remark is
also valid in the case where $\theta= \sigma$). Here, we will prove
the needed results on the Deformed GOE which, thanks to the analysis
of \cite{FePe} and \cite{Ru}, allow us to claim the following
universality.

\begin{theorem}\label{ThmFluctuationsUniversalite} Let $A_N$ be a full
perturbation given by
$(A_N)_{ij}={\theta}/{N}$ for all $(i,j)$. Assume that $\theta>
\sigma$. Let $W_N$ be an arbitrary real Wigner matrix with the
underlying measure $\mu$ being symmetric with a variance $\sigma^2$
and such that $\int|x|^q \, d\mu(x)< + \infty$ for any $q$ in
$\mathbb N$.

Then the largest eigenvalue of the Deformed model $M_N$ has the
Gaussian fluctuations (\ref{FluctDefGOE}).
\end{theorem}

\begin{remark} To be complete, let us notice that the previous result still
holds when we allow the distribution $\nu$ of the diagonal entries of
$W_N$ to be different from~$\mu$ provided that $\nu$ is symmetric
and has moments of any order.
\end{remark}

The fundamental tool of this paper is the Stieltjes transform.
For $z \in{\mathbb C}\setminus{\mathbb R}$, we denote the resolvent
of the matrix
$M_N$ by
\[
G_N(z) = ( zI_N - M_N)^{-1}
\]
and the Stieltjes transform
of the
expectation of the empirical measure of the eigenvalues of $M_N$ by
\[
g_N(z) =
\mathbb{E}(\mathrm{tr}_N (G_N(z))),
\]
where $\mathrm{tr}_N $ is the normalized trace. We also
denote by
\[
g_{\sigma}(z)=\mathbb{E}\bigl((z-s)^{-1}\bigr)
\]
the Stieltjes
transform\footnote{Note that in some papers to which we make
reference, the Stieltjes
transform is defined with the opposite sign.} of a variable $s$ with
semicircular distribution $\mu_{sc}$.

Theorem \ref{ThmASCV} is the analogue of the main statement of
\cite{BS3} established in the context of general spiked population
models. The conclusion of \cite{BS3} requires numerous results
obtained previously by Silverstein and co-authors in \cite{CS},
\cite{BS1} and~\cite{BS2} (a summary of all this literature can be
found in \cite{Bai}, pages 671--675). From very clever and tedious
manipulations of some Stieltjes transforms and the use of the
matricial representation (\ref{spike}), these works highlight a very
close link between the spectra of the Wishart matrices and the
covariance matrix (for quite general covariance matrix which includes
the spiked population model). Our approach mimics the one of
\cite{BS3}. Thus, using the fact that the Deformed Wigner model is
the additive analogue of the spiked population model, several
arguments can be quite easily adapted here (this point has been
explained in Chapter 4 of the Ph.D. thesis \cite{Fe}). Actually, the
main point in the proof consists in establishing that for any
$\varepsilon>0$, almost surely,
%
\begin{equation}\label{inclusion}
\operatorname{Spect}(M_N) \subset K^{\varepsilon}_{\sigma}(\theta_1,
\ldots,\theta_J)
\end{equation}
for all $N$ large, where we have defined
\[
K^{\varepsilon}_{\sigma}(\theta_1, \ldots,\theta_J)=K_{\sigma
}(\theta
_1, \ldots,\theta_J)+(-\varepsilon,\varepsilon)
\]
and
\[
K_{\sigma}(\theta_1, \ldots,\theta_J):=\{\rho_{\theta_J} ;
\ldots; \rho_{\theta
_{J-J_{-\sigma}+1}}\} \cup
[-2 \sigma, 2 \sigma]\cup\{\rho_{\theta
_{J_{+\sigma}}}; \ldots;
\rho_{\theta_1} \}.
\]
This point is the analogue of the main
result of \cite{BS1}. The analysis of \cite{BS1} is based on
technical and numerous considerations of Stieltjes transforms
strongly related to the Wishart context and that cannot be directly
transposed here. Our approach to prove such an inclusion of the
spectrum of $M_N$ is very different from the one of~\cite{BS1}.
Indeed, we use the methods developed by Haagerup and Thorbj\o
rnsen in~\cite{HT}, by Schultz \cite{S} and by the two first
authors of the present article \cite{CD}. The key point of this
approach is to obtain a precise estimation at any point $z \in
\mathbb{C} \setminus\mathbb{R}$ of the following type:
%
\begin{equation}\label{e}
g_{\sigma}(z)- g_N(z)
+ \frac{1}{N}L_{\sigma}(z)= O\biggl(\frac{1}{N^2}\biggr),
\end{equation}
where
$L_{\sigma}$ is the Stieltjes transform of a distribution
$\Lambda_{\sigma}$ with compact support in $K_{\sigma}(\theta_1,
\ldots,\theta_J)$. Indeed such an estimation allows us through the
inverse Stieltjes transform and some variance estimates to deduce
that a.s.,\break $\mathrm{tr}_N (1_{^c K^{\varepsilon}_{\sigma}(\theta_1,
\ldots,\theta_J)}(M_N))=O(N^{-{4}/{3}}) $. Thus the
number of eigenvalues of $M_N$ in $^c
K^{\varepsilon}_{\sigma}(\theta_1, \ldots,\theta_J)$ is almost
surely an $O(N^{-{1}/{3}})$ and since for each $N$ this number
has to be an integer, we deduce
that it is actually equal to zero as $N$ goes to infinity.

Dealing with the particular diagonal perturbation $A_N =
\operatorname{diag}(\theta, 0,\ldots,0)$ such that $\theta>\sigma$, we
obtain the fluctuations of the largest eigenvalue $\lambda_1(M_N)$
(Theorem~\ref{ThmFluctuations}) by an approach close to the one of
\cite{Pa} and the ideas of \cite{BBPbis}. The reasoning relies on
the writing of the rescaled variable $\sqrt N (\lambda_1(M_N)- \rho
_{\theta})$ in terms of the resolvent of a
non-Deformed Wigner matrix. Then, to complete the analysis of \cite
{FePe} and justify Theorem \ref{ThmFluctuationsUniversalite}, we focus
on the particular Deformed GOE model and improve the previous
convergence at the level of Laplace transform.

The paper is organized as follows. In Section \ref{base}, we introduce
preliminary lemmas which will be of basic use later on. Section \ref{section-as} is
devoted to the proof of Theorem \ref{ThmASCV}. We first establish an
equation (called master equation or master inequality) satisfied by
$g_N$ up to some correction of order $\frac{1}{N^2}$ (see\vspace*{2pt} Section
\ref{TheMasterequation}). Then we explain how this master equation gives rise to an
estimation of type (\ref{e}) and thus to the inclusion
(\ref{inclusion}) of the spectrum of $M_N$ in
$K^{\varepsilon}_{\sigma}(\theta_1, \ldots,\theta_J)$ (see
Sections~\ref{ssu} and \ref{incsp}). In Section \ref{asc}, we use this inclusion to relate the
asymptotic spectra of~$A_N$ and $M_N$ and then deduce Theorem
\ref{ThmASCV}. Section \ref{sec5} deals with the fluctuations results. The
proof of Theorem \ref{ThmFluctuations} is given in Section
\ref{sectionThmfluctu};
Theorem~\ref{ThmFluctuationsUniversalite} is justified in Section \ref{sectiongoe}.

\section{Basic lemmas}\label{base}

We assume that the distribution $\mu$ of the entries of the
Wigner matrix $W_N$ satisfies a Poincar\'e
inequality: there exists a positive constant~$C$ such that for any
$\mathcal{C}^1$ function $f\dvtx {\mathbb R}\rightarrow{\mathbb C}$ such that
$f$ and
$f'$ are in $L^2(\mu)$,
\[
\mathbf{V}(f)\leq C\int\vert f' \vert^2 \,d \mu,
\]
with $\mathbf{V}(f) = \mathbb{E}(\vert
f-\mathbb{E}(f)\vert^2)$.

Let $\mathrm{Tr}$ denote the classical trace.

For any matrix $M$, define $\Vert M\Vert_2 = (\operatorname{Tr}(M^*M))^{1/2}$
the Hilbert--Schmidt norm. Let $\Psi\dvtx (M_N({\mathbb C})_{sa})
\rightarrow
\mathbb{R}^{N^2}$ [resp., $\Psi\dvtx (M_N({\mathbb C})_{s}) \rightarrow
\mathbb{R}^{{N(N+1)}/{2}}$] be the canonical isomorphism which
maps an Hermitian (resp., symmetric) matrix $M$ to the real parts and
the imaginary parts
of its entries (resp., to the entries)
$(M)_{ij}, i\leq j$.

\begin{lemma}\label{variance} Let $M_N$ be the complex (resp., real)
Wigner Deformed matrix introduced in Section \ref{results}.
For any $\mathcal{C}^1$ function $f\dvtx {\mathbb R}^{N^2}$ (resp., $\mathbb
{R}^{{N(N+1)}/{2}} ) \rightarrow{\mathbb C}$ such that $f$ and the gradient
$\nabla(f )$ are both polynomially
bounded,
%
\begin{equation}\label{Poincare}\mathbf{V}{[ f\circ\Psi(M_N
)]} \leq\frac{C}{N}\mathbb{E}\{\Vert\nabla[f
\circ\Psi(M_N)] \Vert_2^2
\}.
\end{equation}
\end{lemma}

\begin{pf}
According to Lemma 3.2 in \cite{CD},
%
\begin{equation}\label{PoincareX}\mathbf{V}{[ f\circ\Psi(X_N
)]} \leq\frac{C}{N}\mathbb{E}\{\Vert\nabla[f
\circ\Psi(X_N)] \Vert_2^2
\}.
\end{equation}
Note that even if the result in \cite{CD} is stated in the Hermitian
case, the proof is valid and the result still holds in the symmetric case.
Now (\ref{Poincare}) follows putting $g(x_{ij};i\leq j ):=f(x_{ij}+(A_N)_{ij}
;i\leq j)$ in (\ref{PoincareX}) and noticing that the $(A_N)_{ij}$ are
uniformly bounded in $i,j,N$.
\end{pf}

This lemma will be useful to estimate many variances. Now,
we recall some useful properties of the resolvent (see \cite{KKP,CD}).

\begin{lemma} \label{lem0}
For an $N \times N$ Hermitian or symmetric matrix $M$, for any $z \in
{\mathbb C} \setminus\operatorname{Spect}(M)$, we denote by $G(z) :=
(zI_N-M)^{-1}$ the resolvent of $M$.

Let $z \in{\mathbb C} \setminus{\mathbb R}$.
\begin{longlist}[(iii)]
\item[(i)]  $\Vert G(z) \Vert\leq|\Im m(z)|^{-1}$ where $\Vert\cdot
\Vert
$ denotes the operator norm.
\item[(ii)]  $\vert G(z)_{ij} \vert\leq|\Im m(z)|^{-1}$ for all $i,j =
1, \ldots, N$.
\item[(iii)]  For $p\geq2$,
\[
\frac{1}{N} \sum_{i,j = 1}^N \vert G(z)_{ij} \vert^p \leq(|\Im
m(z)|^{-1})^p.
\]
\item[(iv)]  The derivative with respect to $M$ of the resolvent $G(z)$
satisfies
\[
G'_M(z)\cdot B = G(z)BG(z)  \qquad\mbox{for any matrix $B$}.
\]
\item[(v)]
Let $z \in
{\mathbb C}$ such that
$|z| > \Vert M \Vert$; we have
\[
\Vert G(z) \Vert\leq\frac{1}{|z| - \Vert M \Vert}.
\]
\end{longlist}
\end{lemma}

\begin{pf}
We just mention that (v) comes readily noticing that the
eigenvalues of the normal matrix $G(z)$ are the $ \frac{1}{z -
\lambda_i( M )}$, $i=1, \ldots, N.$
\end{pf}

We will also need the following estimations on the
Stieltjes transform $g_\sigma$ of the semicircular distribution
$\mu_{sc}$.

\begin{lemma}\label{lem0bis}
$g _{\sigma}$ is analytic on $\mathbb C \setminus[-2 \sigma,2
\sigma]$ and
\begin{longlist}
\item $\forall z \in\{ z \in\mathbb C\dvtx  \Im m (z) \not= 0 \}$,
%
\begin{eqnarray}
\label{g}
\sigma^2 g_\sigma^2(z) -z g_\sigma(z) + 1 &=& 0,  \\
\label{bornesupg}
\vert g _{\sigma}(z) \vert
&\leq&
\vert\Im m (z) \vert^{-1},  \\
\label{gmoins}
\vert g_\sigma(z)^{-1}\vert
&\leq&|z|+ \sigma^2 \vert\Im m (z) \vert
^{-1}, \\
\label{bgprime}
\vert g
'_{\sigma}(z)\vert
&=& \biggl\vert\int\frac{1}{(z-t)^2} \, d\mu_{\sigma} (t)
\biggr\vert\leq|\Im m(z)|^{-2},
\\[-3pt]
\label{0.1}
&&\hspace*{-86.3pt}\mbox{for $a>0$, $\theta\in{\mathbb R}$,}\qquad
\biggl\vert\frac{1}{ag_{\sigma}(z) -z + \theta} \biggr\vert\leq
|\Im m(z)|^{-1}.
\end{eqnarray}
\item $\forall z \in\{ z \in\mathbb C\dvtx   \vert z \vert>
2\sigma\}$,
%
\begin{eqnarray}
\label{bgout}
\vert g _{\sigma}(z) \vert
&\leq&
\frac{1}{|z| - 2 \sigma}.  \\
\label{bgprimeout}
\vert g
'_{\sigma}(z)\vert&=& \biggl\vert\int\frac{1}{(z-t)^2} \, d\mu_{\sigma} (t)
\biggr\vert\leq
\frac{1}{(|z| - 2 \sigma)^2},  \\
\label{bgmoinsout}
\vert g _{\sigma}(z) \vert^{-1} &\leq&|z|+ \frac{\sigma^2}{|z| - 2
\sigma}.
\end{eqnarray}
\end{longlist}
\end{lemma}
\begin{pf}
For (\ref{g}), we refer the reader to Section 3.1 of
\cite{Bai}. Equation (\ref{0.1}) is a consequence of $\Im
m(g_{\sigma}(z))\Im m(z) <0$. Other
inequalities derive from \textup{(\ref{g})} and the definition of
$g_{\sigma}.$
\end{pf}

\section{Almost sure convergence of the first extremal
eigenvalues}\label{section-as}

Sections \ref{TheMasterequation}, \ref{ssu} and \ref{incsp}
below describe the different steps of the proof of the inclusion
(\ref{inclusion}). We choose to develop the case of the complex
Deformed Wigner model and just to point out some differences with the
real model case (at the end of Section \ref{incsp}) since the approach
would be basically the same. In these sections,
we will often refer the reader to the paper \cite{CD} where the authors
deal with several independent non-Deformed Wigner matrices. The reader
needs to fix $r=1$, $m=1$, $a_0=0$, $a_1 =\sigma$ and to change the
notation $\lambda=z$, $G_N=g_N$, $G=g_{\sigma}$ in \cite{CD} in order
to use the different proofs we refer to in the present framework.
We shall denote by $P_k$ any polynomial of degree $k$ with positive
coefficients and by $C$, $K$ any constants; $P_k$, $C$, $K$ can depend
on the fixed eigenvalues of $A_N$ and may vary from line to line. We
also adopt the following convention to simplify the writing: we
sometimes state in the proofs below that a quantity $\Delta_N(z)$, $z
\in{\mathbb C}\setminus{\mathbb R}$ is $O(N^{-p})$, $p=1,2$. This
means precisely that
\[
|\Delta_N(z)| \leq(\vert z \vert+K)^l \frac{P_k(|\Im m(z)|^{-1})}{N^p}
\]
for some $k$ and some $l$ and we give the precise majoration in the
statements of the theorems or propositions.

Section \ref{asc} explains how to deduce Theorem \ref{ThmASCV} from the
inclusion (\ref{inclusion}).

The goal of Sections \ref{TheMasterequation} and \ref{ssu} is to
establish Proposition \ref{propestimdifL2} below which is fundamental
in the proof of the inclusion
(\ref{inclusion}).
Before describing rigorously the different ideas of these two
sections, let us
help the reader's intuition by a heuristic understanding of the approach.
Assume that we can establish that $g_N(z)$ satisfied the rough
quadratic equation (also called master inequality):
\[
\sigma^2 g_N^2(z) -z g_N(z) + 1 + \frac{1}{N}
E_\sigma(z)
= O\biggl(\frac{1}{N^2}\biggr).
\]
Then, for any suitable $z$, divided by $g_N(z)$ the last approximation
would provide us
an estimation of
$\Lambda_N(z) -z$ where $\Lambda_N(z)= z_{\sigma}(g_N(z))$ with
$z _{\sigma}(g)= \frac{1}{g} + \sigma^2 g$ being the inverse function
of $g _{\sigma}$
(see Lemma \ref{lemme-g-sigma} below).
Then, intuitively, a Taylor expansion of $g _{\sigma}$ between
$\Lambda
_N(z)$ and $z$ would lead to an estimation of
the type
\[
g_N(z) -g _{\sigma}(z) = - \frac{1}{N}
(g_{\sigma}(z))^{-1}g'_{\sigma}(z) E_\sigma(z)+ O\biggl( \frac{1}{N^2}\biggr).
\]
This intuitive process may throw light on the expression (\ref
{Lsigma1}) of
$L_{\sigma}(z)$ in Proposition \ref{propestimdifL2} below.
\eject

\subsection{The master equation}\label{TheMasterequation}
\subsubsection{A first master inequality}
In order to obtain a master equation for~$g_N(z)$, we first consider
the Gaussian case, that is, $X_N=X_N^G$ is distributed as the $\operatorname{GUE}(N,
\sigma
^2/N$) distribution.\footnote{Throughout this section, we will drop the
subscript $G$ in the interest of clarity.}

Let us recall the integration by parts formula for the Gaussian
distribution.
\begin{lemma} \label{lemIPP}
Let $\Phi$ be a complex-valued $\mathcal C^1$ function on
$(M_N({\mathbb C})_{sa})$ and $X_N \sim \operatorname{GUE}(N, {\sigma^2}/{N}$). Then,
%
\begin{equation} \label{IPP1}
\mathbb{E}[\phi'(X_N)\cdot H] = \frac{N}{\sigma^2}\mathbb{E}[ \phi
(X_N) \operatorname{Tr}
(X_N H)]
\end{equation}
for any Hermitian matrix $H$, or by linearity for $H = E_{jk}$, $1\leq
j, k \leq N$ where $E_{jk}$, $1\leq j, k\leq N$ is the canonical basis
of the complex space of $N \times N$ matrices.
\end{lemma}

We\vspace*{1pt} apply the above lemma to the function $\Phi(X_N) = (G_N(z))_{ij} =
((zI_N - X_N -A_N)^{-1})_{ij}$, $z
\in{\mathbb C} \setminus{\mathbb R}$, $1\leq i, j \leq N$. In order
to simplify the notation, we write
$(G_N(z))_{ij} = G_{ij}$. We obtain, for $H= E_{ij}$,
\begin{eqnarray*}
\mathbb{E}((GHG)_{ij})& = &\frac{N}{\sigma^2}\mathbb{E}[G_{ij}
\operatorname{Tr}(X_N
H)], \\
\mathbb{E}(G_{ii} G_{jj}) &=& \frac{N}{\sigma^2} \mathbb{E}[G_{ij}
(X_N )_{ji}].
\end{eqnarray*}
Now, we consider the normalized sum $\frac{1}{N^2} \sum_{ij}$ of the
previous identities to obtain
\[
\mathbb{E}((\mathrm{tr}_N G)^2) = \frac{1}{\sigma^2}
\mathbb{E}(\mathrm{tr}_N(GX_N)).
\]
Then, since
\[
GX_N = (z-X_N - A_N)^{-1}(X_N + A_N -zI_N -A_N +zI_N) = -I_N -GA_N + zG,
\]
we obtain the following master equation:
\[
\mathbb{E}((\mathrm{tr}_N G)^2) + \frac{1}{\sigma^2} \bigl(- z \mathbb
{E}(\mathrm{tr}_N G) +
1 + \mathbb{E}(\mathrm{tr}_N GA_N)\bigr) = 0.
\]
Now, it is well known (see \cite{CD,HT} and Lemma \ref
{variance}) that
\[
\operatorname{Var}(\mathrm{tr}_N(G)) \leq\frac{C|\Im m (z)|^{-4}}{N^2}.
\]
Thus, in the case where $X_N=X_N^G$ we obtain:
\begin{proposition}
The Stieltjes transform $g_N$ satisfies the following inequality:
%
\begin{equation} \label{masteqGUE}
\biggl| \sigma^2 g_N^2(z) -z g_N(z) + 1 + \frac{1}{N} \mathbb
{E}(\operatorname{Tr}(G_N(z)
A_N))\biggr| \leq C \frac{|\Im m (z)|^{-4}}{N^2}.
\end{equation}
Note that since $A_N$ is of finite rank, $\mathbb{E}(\operatorname
{Tr}(G_N(z)A_N))
\leq C$ where $C$ is a constant independent of $N$ (depending on the
eigenvalues of $A_N$ and $z$).
\end{proposition}

We now explain how to obtain the corresponding \textup{(\ref
{masteqGUE})} in the Wigner case. Since the computations are the same as
in \cite{CD}\footnote{This paper treats the case of several independent
non-Deformed Wigner matrices.} and \cite{KKP},\footnote{The authors
considered a non-Deformed Wigner matrix in the symmetric real
setting.} we just give some hints of the proof.

\begin{step}\label{step1}
The integration by parts formula for the Gaussian
distribution is replaced by the following tool:
\begin{lemma} \label{lem1}
Let $\xi$ be a real-valued random variable such that $\mathbb
{E}(\vert
\xi
\vert^{p+2})<\infty$. Let $\phi$ be a function from ${\mathbb R}$ to
${\mathbb C}$
such that the first $p+1$ derivatives are continuous and bounded.
Then,
%
\begin{equation}\label{IPP2}\mathbb{E}(\xi\phi(\xi)) = \sum_{a=0}^p
\frac{\kappa_{a+1}}{a!}\mathbb{E}\bigl(\phi^{(a)}(\xi)\bigr) +
\epsilon
\end{equation}
where $\kappa_{a}$ are the cumulants of
$\xi$, $\vert\epsilon\vert\leq C \sup_t \vert
\phi^{(p+1)}(t)\vert\mathbb{E}(\vert\xi\vert^{p+2})$, $C$
depends on~$p$ only.
\end{lemma}

We apply this lemma with the function $\phi(\xi)$ given, as before,
by $\phi(\xi) = G_{ij}$ and~$\xi$ is now one of the variables $\Re
e((X_N)_{kl})$, $\Im m((X_N)_{kl})$. Note that, since the above
random variables are symmetric, only the odd derivatives in
\textup{(\ref{IPP2})} give a nonnull term. Moreover, as we are
concerned by
estimation of order $\frac{1}{N^2}$ of $g_N$, we only need to
consider \textup{(\ref{IPP2})} up to the third derivative (see
\cite{CD}).\vadjust{\goodbreak}
The computation of the first derivative will provide the same term
as in the Gaussian case.\looseness=1
\end{step}

\begin{step}
Study of the third derivative.

We refer to \cite{CD} or \cite{KKP} for a detailed study of the third
derivative. Using some bounds on $G_N$ (see Lemma \ref{lem0}), we can
prove that the only term arising from the third derivative in the
master equation, giving a contribution of order $\frac{1}{N}$, is
\[
\frac{1}{N} \mathbb{E}\Biggl[ \Biggl( \frac{1}{N} \sum_{k=1}^N
G_{kk}^2\Biggr)^2\Biggr].
\]
In conclusion, the first master equation in the Wigner case reads as follows:
\begin{theorem} \label{theomasteq1}
For $ z \in{\mathbb C} \setminus{\mathbb R}$,
$g_N(z)$ satisfies
%
\begin{eqnarray} \label{masteq}
&&\Biggl| \sigma^2 g_N(z)^2 -z g_N(z)
+ 1 +  \frac{1}{N} \mathbb{E}[\operatorname
{Tr}(G_N(z)A_N)]\nonumber\\[-8pt]\\[-8pt]
&&\qquad\hspace*{43.1pt}
+ \frac{1}{N} \frac{\kappa_4}{2}
\mathbb{E}\Biggl[ \Biggl( \frac{1}{N} \sum_{k=1}^N
(G_N(z))_{kk}^2
\Biggr)^2\Biggr] \Biggr|\leq\frac{P_6(|\Im m (z)|^{-1})}{N^2},\nonumber
\end{eqnarray}
where $\kappa_4$ is the fourth cumulant of the distribution $\mu$.
\end{theorem}
\end{step}

\subsubsection{Estimation of $|g_N -g_{\sigma}|$}
Since
\[
|\mathbb{E}[\operatorname{Tr}(G_N(z)A_N)]| \vee\Biggl|\mathbb{E}\Biggl[
\Biggl( \frac{1}{N}
\sum_{k=1}^N (G_N(z))_{kk}^2\Biggr)^2\Biggr] \Biggr| \leq P_4( |\Im m (z)|^{-1}),
\]
Theorem \ref{theomasteq1} implies that for any $z \in{\mathbb
C}\setminus
{\mathbb R}$,
%
\begin{equation} \label{estenunsurn}
|\sigma^2 g_N(z)^2 - z g_N(z) +1| \leq\frac{P_6(|\Im m
(z)|^{-1})}{N}.
\end{equation}
To estimate $|g_N -g_{\sigma}|$ from (\ref{g}) and (\ref{estenunsurn}), we follow the method initiated in
\cite{HT} and~\cite{S}. We do not develop it here since it follows
exactly the lines of Section 3.4 in \cite{CD}
but we briefly recall the main arguments and results which will be
useful later on.
We define the open connected set
\[
\mathcal{O}'_N= \biggl\{z \in{\mathbb C},   \Im m (z) >0, \frac
{P_6(|\Im m
(z)|^{-1})}{N}(\sigma^2 |\Im m (z)|^{-1}+| z|)< \frac{1}{4|\Im m
(z)|^{-1}}\biggr\}.
\]
For any $z$ in ${\mathbb C}$ such that $\Im m (z) >0,$ we set
%
\begin{equation} \label{deflambdaN}
\Lambda_N(z):= \sigma^2 g_N(z) +\frac{1}{g_N(z)}.
\end{equation}
One can prove that for any $z$ in $\mathcal{O}'_N$:
\begin{itemize}
\item$g_N(z)\neq0$ and
%
\begin{equation} \label{gNmoins}\frac{1}{|g_N(z)|} \leq2 (\sigma^2
|\Im m (z)|^{-1}+| z|);
\end{equation}
\item from (\ref{estenunsurn}) and (\ref{gNmoins}),
%
\begin{equation} \label{lambdaN}
|\Lambda_N(z)-z|\leq\frac{P_6(|\Im m (z)|^{-1})}{N}2\bigl(\sigma^2 |\Im m
(z)|^{-1}+| z|\bigr)
\end{equation}
and
%
\begin{equation} \label{binfLambda} \Im m (\Lambda_N(z))\geq\frac
{\Im
m (z)}{2}>0;
\end{equation}
\item writing (\ref{g}) at the point $\Lambda_N(z)$, we
easily get that
%
\begin{equation} \label{id}g_N(z)=g_\sigma(\Lambda_N(z))
\end{equation}
on the nonempty open subset
$\mathcal{O}''_N=\{z \in\mathcal{O}'_N,\Im m (z) > \sqrt{2}
\sigma\}$
and then on $\mathcal{O}'_N$ by the principle of uniqueness of continuation.
\end{itemize}
Using
\begin{eqnarray*}
\vert g_N(z)-g_\sigma(z)\vert& = & \bigl\vert\mathbb{E}\bigl[
(z-s)^{-1}\bigl(\Lambda_N(z)-s\bigr)^{-1}\bigl(\Lambda_N(z)-z\bigr) \bigr] \bigr\vert\\
& \leq& \Im m(z)\cdot  \Im m(\Lambda_N(z))\cdot \vert\Lambda_N(z)-z \vert,
\end{eqnarray*}
this
allows us to get an estimation of $|g_N(z)-g_\sigma(z)|$ on $\mathcal
{O}'_N$ and then to deduce:
\begin{proposition}
For any $z \in{\mathbb C}$ such that $\Im m(z)>0$,
%
\begin{equation}\label{g-gN}
|g_N(z) - g_\sigma(z)| \leq(|z|+K)\frac{P_9(|\Im m (z)|^{-1})}{N}.
\end{equation}
\end{proposition}

\subsubsection{Study of the additional term $\mathbb{E}[\operatorname
{Tr}(A_N G_N(z))]$}
From now on and until the end of Section \ref{TheMasterequation}, we
denote by $\gamma_1, \ldots, \gamma_r$ the nonnull eigenvalues of
$A_N$ ($\gamma_i = \theta_j$ for some $j \not= j_0$) in order to
simplify the writing. Let $U_N:= U$ be a unitary matrix such that
$A_N = U^* \Delta U$ where $\Delta$ is the diagonal matrix with
entries $\Delta_{ii} = \gamma_i, \  i\leq r;   \Delta_{ii} = 0,\
i > r$. We set
%
\begin{equation} \label{defhN}
h_N(z) := \mathbb{E}[\operatorname{Tr}(A_N G_N(z))] = \sum_{k=1}^r
\gamma_k
\sum_{i,j = 1}^N U^*_{ik} U_{kj} \mathbb{E}[G_{ji}].
\end{equation}
Our aim is to express $h_N(z)$ in terms of the Stieltjes transform
$g_N(z)$ for $N$ large, using the integration by parts formula. Note
that since we want an estimation of order $O(N^{-2})$ in the
master inequality \textup{(\ref{masteq})}, we only need an estimation of
$h_N(z)$ of order $O(N^{-1})$. As in the previous subsection,
we first write the equation in the Gaussian case and then study the
additional term (third derivative) in the Wigner case.\vspace*{12pt}

(a) \textit{Gaussian case}.
Apply \textup{(\ref{IPP1})} to $\Phi(X_N) = G_{jl}$ and
$H = E_{il}$
to get
\[
\mathbb{E}[G_{ji} G_{ll}] = \frac{N}{\sigma^2} \mathbb{E}[G_{jl}(X_N)_{li}]
\]
and
\[
\frac{1}{N} \sum_{l=1}^N \mathbb{E}[G_{ji} G_{ll}] = \frac
{1}{\sigma^2}
\mathbb{E}[(GX_N)_{ji}].
\]
Expressing $GX_N$ in terms of $GA_N$, we obtain
%
\begin{equation} \label{ijfix}
I_{ji} := \sigma^2 \mathbb{E}[ G_{ji} \mathrm{tr}_N(G)] + \delta_{ij} -
z\mathbb{E}[G_{ji}] + \mathbb{E}[(GA_N)_{ji}] = 0.
\end{equation}
Now, we consider the sum $\sum_{i,j} U^*_{ik}U_{kj} I_{ji}$, $k =1,
\ldots, r$ fixed and we denote
$\alpha_k = \sum_{i,j} U^*_{ik}U_{kj} G_{ji} = (UGU^*)_{kk}$.
Then, we have the following equality, using that $U$ is unitary:
\[
\sigma^2 \mathbb{E}[\alpha_k \mathrm{tr}_N(G)] + 1 - z \mathbb
{E}[\alpha_k] +
\sum_{i,j} U^*_{ik}U_{kj} \mathbb{E}[(GA_N)_{ji}] = 0.
\]
Now,
\begin{eqnarray*}
\sum_{i,j} U^*_{ik}U_{kj} \mathbb{E}[(GA_N)_{ji}] &= &\mathbb
{E}[(UGA_NU^*)_{kk}]\\[-5pt]
&=& \mathbb{E}[(UGU^*\Delta UU^*)_{kk}] \\
&=&
\gamma_k \mathbb{E}[(UGU^*)_{kk}] = \gamma_k \mathbb{E}[\alpha_k].
\end{eqnarray*}
Therefore,
\[
\sigma^2 \mathbb{E}[\alpha_k \mathrm{tr}_N(G)] + 1 + ( \gamma_k -z)
\mathbb
{E}[\alpha_k] = 0.
\]
Since $\alpha_k$ is bounded and $\operatorname{Var}(\mathrm{tr}_N(G)) =
O(N^{-2})$, we obtain
%
\begin{equation} \label{eqalpha}
\mathbb{E}[\alpha_k]\bigl( \sigma^2 g_N(z) + \gamma_k -z\bigr) + 1 = O\biggl(\frac{1}{N}\biggr).
\end{equation}
Then using (\ref{g-gN}) we deduce that
$\mathbb{E}[\alpha_k]( \sigma^2 g_\sigma(z) + \gamma_k -z) + 1 =
O(N^{-1}) $ and using (\ref{0.1}):
%
\begin{equation} \label{eqhN}
h_N(z) = \sum_{k=1}^r \gamma_k \mathbb{E}[\alpha_k] = \sum_{k=1}^r
\frac{
\gamma_k}{z - \sigma^2 g _{\sigma}(z) - \gamma_k} +
O\biggl(\frac{1}{N}\biggr).\vadjust{\goodbreak}
\end{equation}

(b) \textit{The general Wigner case.}
We shall prove that \textup{(\ref{eqalpha})} still holds. We now rely on
Lemma \ref{lem1} to obtain the analogue of \textup{(\ref{ijfix})}:
%
\begin{eqnarray} \label{ijfix1}
J_{ij}:\!&=& \sigma^2 \mathbb{E}[ G_{ji} \mathrm{tr}_NG] + \delta_{ij} -
z\mathbb{E}[G_{ji}]\nonumber\\
&&{} + \mathbb{E}[(GA_N)_{ji}] +
\frac{\kappa_4}{6N^2} \sum_{l=1}^{N}
\mathbb{E}[A_{i,j,l}]\\
&=&
O\biggl(\frac{1}{N^2}\biggr).\nonumber
\end{eqnarray}
The term $A_{i,j,l}$ is a fixed linear combination of the third
derivative of $\Phi:= G_{jl}$ with respect to
$\mathcal{R}e (X_N)_{il}$ (i.e., in the direction $e_{il} = E_{il} + E_{li})$ and
$\Im m(X_N)_{il}$ [i.e., in the direction $f_{il}:= \sqrt{-1} (E_{il}
-E_{li})$]. We do not need to write the exact form of this term since we
just want to show that this term will give a contribution of order
$O(N^{-1})$ in the equation for $h_N(z)$.
Let us write the derivative in the direction $e_{il}$:
\[
\mathbb{E}[(Ge_{il}Ge_{il}Ge_{il}G)_{jl}]
\]
which is the sum of eight terms of the form
%
\begin{equation} \label{derivee3}
\mathbb{E}[ G_{ji_1} G_{i_2 i_3} G_{i_4 i_5} G_{i_6 l} ],
\end{equation}
where if $i_{2q+1} = i$ (resp., $l$), then $i_{2q+2} = l$ (resp., $i$),
$q = 0,1,2$.
\begin{lemma} \label{lem2} Let $1 \leq k \leq r$ fixed; then
%
\begin{equation}
F(N) :=\Biggl| \sum_{i,j =1}^N U^*_{ik} U_{kj} \frac{1}{N} \sum_{l=1}^N
\mathbb{E}[A_{i,j,l}] \Biggr| \leq C |\Im m (z)|^{-4}
\end{equation}
for a numerical constant $C$.
\end{lemma}

\begin{pf}
$F(N)$ is the sum of eight terms corresponding to \textup{(\ref
{derivee3})}. Let
us write, for example, the term corresponding to $i_1= i$, $i_3 = i$,
$i_5 =i$:
\begin{eqnarray*}
\frac{1}{N} \sum_{i,j,l} U^*_{ik} U_{kj} E[ G_{ji}G_{li}G_{li}G_{ll}]
& = & \mathbb{E}\Biggl[ \frac{1}{N} \sum_{i,l} U^*_{ik}
(UG)_{ki}G_{li}G_{li}G_{ll}\Biggr]\\
& = & \mathbb{E}\Biggl[ \frac{1}{N} \sum_{i} U^*_{ik}
(UG)_{ki}(G^TG^DG^T)_{ii}\Biggr],
\end{eqnarray*}
where the superscript $T$ denotes the transpose of the matrix and $G^D$
is the diagonal matrix with entries $G_{ii}$. From the bounds $\Vert
G_N(z) \Vert\leq|\Im m (z)|^{-1}$ and $\Vert U\Vert= 1$, we get the
bound given in the lemma.\vadjust{\goodbreak}

We give the majoration for the term corresponding to $i_1= l$, $i_3 =
l$, $i_5 =l$:
\begin{eqnarray*}
\frac{1}{N} \sum_{i,j,l} U^*_{ik} U_{kj} \mathbb{E}[ G_{jl}G_{il}^3] =
\mathbb{E}\Biggl[ \frac{1}{N} \sum_{i,l} U^*_{ik}
(UG)_{kl}G_{il}^3\Biggr].
\end{eqnarray*}
Its absolute value is bounded by $\mathbb{E}[ \frac{1}{N} \sum
_{i,l} |G_{il}| ^3] |\Im m (z)|^{-1}$ and thanks to Lemma~\ref
{lem0} by $|\Im m (z)|^{-4}$. The other terms are treated in the same
way.
\end{pf}

As in the Gaussian case, we now consider the sum $\sum_{i,j}
U^*_{ik}U_{kj} J_{ji}$. From Lemma \ref{lem2} and the bound (using the
Cauchy--Schwarz inequality)
\[
\sum_{i,j =1}^N |U^*_{ik} U_{kj}| \leq N,
\]
we still get \textup{(\ref{eqalpha})} and thus \textup{(\ref
{eqhN})}. More precisely, we proved:
\begin{proposition} For any $z \in{\mathbb C}$ such that $\Im m(z)>0$,
\[
\Biggl| \mathbb{E}[\operatorname{Tr}(A_NG_N(z))] - \sum_{k=1}^r \frac
{\gamma_k}{ z -
\sigma
^2 g_\sigma(z) - \gamma_k} \Biggr| \leq\frac{P_{11}(|\Im m (z)|^{-1})}{N} (K
+ |z|).
\]
\end{proposition}

\subsubsection{Convergence of $\mathbb{E}[(\frac{1}{N}\sum_{k=1}^N G_{kk}^2)^2]$}

We now study the last term in the master inequality of Theorem \ref
{theomasteq1}. For the non-Deformed
Wigner matrices, it is shown in \cite{KKP} that
\[
R_N(z) := \mathbb{E}\Biggl[ \Biggl( \frac{1}{N} \sum_{k=1}^N
G_{kk}^2\Biggr)^2\Biggr] \mathop{\longrightarrow}_{N \mathop
{\longrightarrow}\infty} g_\sigma^4(z).
\]
Moreover, Proposition 3.2 in \cite{CD}, in the more general setting of
several independent Wigner matrices, gives
an estimate of $ | R_N(z) - g_\sigma^4(z) | $.
The above convergence holds true in the Deformed case. We just give
some hints of the proof of the estimate of $ | R_N(z) - g_\sigma^4(z) |
$ since the computations are almost the same as in the non-Deformed
case. Let us set
\[
d_N(z) = \frac{1}{N} \sum_{k=1}^N G_{kk}^2.
\]
We start from the resolvent identity
\begin{eqnarray*}
z G_{kk} &= & 1 + \sum_{l=1}^N (M_N)_{kl} G_{lk} \\
&=& 1 + \sum_{l=1}^N (A_N)_{kl} G_{lk} + \sum_{l=1}^N (X_N)_{kl} G_{lk}
\end{eqnarray*}
and
\[
z d_N(z) = \frac{1}{N} \sum_{k=1}^N G_{kk} + \frac{1}{N} \sum_{k=1}^N
(A_N G)_{kk}G_{kk} +
\frac{1}{N} \sum_{k, l=1}^N (X_N)_{kl} G_{lk}G_{kk}.
\]
For the last term, we apply an integration by parts formula (Lemma \ref
{lem1}) to obtain (see \cite{KKP,CD})
\[
\mathbb{E}\Biggl[ \frac{1}{N} \sum_{k, l=1}^N (X_N)_{kl} G_{lk}G_{kk}
\Biggr] = \sigma^2
\mathbb{E}\Biggl[ \Biggl(\frac{1}{N} \sum_{k=1}^N G_{kk} \Biggr)\, d_N(z) \Biggr] +
O\biggl(\frac{1}{N}\biggr).
\]
It remains to see that the additional term due to $A_N$ is of order
$O(N^{-1})$:
\[
\frac{1}{N} \sum_{k=1}^N (A_N G)_{kk}G_{kk} = \frac{1}{N} \sum_{p=1}^r
\gamma_p (UGG^DU^*)_{pp}
\]
and
\[
\Biggl| \frac{1}{N} \sum_{k=1}^N (A_N G)_{kk}G_{kk} \Biggr| \leq\Biggl(\sum_{p=1}^r |
\gamma_p| \Biggr)\frac{|\Im m (z)|^{-2}}{N}.
\]
We thus obtain (again with the help of a variance estimate)
\[
z \mathbb{E}[ d_N(z)] = g_N(z) + \sigma^2 g_N(z) \mathbb{E}[ d_N(z)] +
O\biggl(\frac{1}{N}\biggr).
\]
Then using (\ref{g-gN}) and since $d_N(z)$ is bounded we deduce that
\[
z \mathbb{E}[ d_N(z)] = g_\sigma(z) + \sigma^2 g_\sigma(z) \mathbb{E}[
d_N(z)] + O\biggl(\frac{1}{N}\biggr).
\]
Thus [using (\ref{0.1})]
\[
\mathbb{E}[ d_N(z)] = \frac{g_\sigma(z)}{z - \sigma^2 g_\sigma(z)} +
O\biggl(\frac{1}{N}\biggr) \mathop{\longrightarrow}_{N \mathop{\longrightarrow
}\infty} \frac{g_\sigma(z)}{z -
\sigma^2
g_\sigma(z)} = g_\sigma^2 (z).
\]
Now, using some variance estimate,
\[
\mathbb{E}[d_N^2(z)] = (\mathbb{E}[ d_N(z)])^2 + O\biggl(\frac{1}{N}\biggr) =
g_\sigma^4(z) + O\biggl(\frac{1}{N}\biggr).
\]
We can now give our final master inequality for $g_N(z)$ following our
previous estimates:
\begin{theorem} \label{theomasteq2}
For $ z \in{\mathbb C}$ such that $\Im m(z) >0$, $g_N(z)$ satisfies
\[
\biggl| \sigma^2 g_N^2(z) -z g_N(z) + 1 + \frac{1}{N}
E_\sigma(z)
\biggr| \leq\frac{P_{14}(|\Im m (z)|^{-1})}{N^2}(|z| + K),
\]
where ${E_\sigma(z)=\sum_{k=1}^r \frac{ \gamma_k}{ z -
\sigma^2 g_\sigma(z) - \gamma_k} + \frac{\kappa_4}{2}
g_\sigma^4(z)},$
$\kappa_4$ is the fourth cumulant of the distribution $\mu$.
\end{theorem}

Note that $E_\sigma(z)$ can be written in terms of the distinct
eigenvalues $\theta_j$ of $A_N$ as
%
\begin{equation}\label{Esigma1} E_\sigma(z)=\sum_{j=1, j \not= j_0}^J
k_j \frac{ \theta_j}{ z - \sigma^2 g_\sigma(z) - \theta_j} + \frac
{\kappa_4}{2}
g_\sigma^4(z).
\end{equation}

Let us set
%
\begin{equation}\label{Lsigma1}
L_{\sigma}(z)= g_\sigma(z)^{-1} {\mathbb E}\bigl((z-s)^{-2}\bigr)E_\sigma(z),
\end{equation}
where $s$ is a centered semicircular random variable with variance
$\sigma^2$.

\subsection{Estimation of $\vert g_\sigma(z)-g_N(z)+\frac{1}{N}L_{\sigma}(z)\vert$}\label{ssu}

The method\vspace*{2pt} is roughly the same as the one described
in Section 3.6 in \cite{CD}. Nevertheless
we choose to develop it here for the reader's convenience.
We have for any
$z$ in
$\mathcal{O}'_n$, by using (\ref{deflambdaN}) and (\ref{id}),
\begin{eqnarray*}
&&\biggl\vert g_\sigma(z)-g_N(z) +\frac{1}{N} L_{\sigma}(z) \biggr\vert\\
&&\qquad= \biggl\vert g_\sigma(z)-g_{\sigma}(\Lambda_N(z))
+\frac
{1}{N} L_{\sigma}(z) \biggr\vert\\
&&\qquad=\biggl\vert\mathbb{E}\biggl[ (z-s)^{-1}\bigl(\Lambda_N(z)-s\bigr)^{-1}\bigl(\Lambda
_N(z)-z\bigr) + \frac{1}{N}g_\sigma(z)^{-1} (z-s)^{-2}E_\sigma(z)
\biggr]\biggr\vert
\\
&&\qquad\leq\biggl\vert\mathbb{E} \biggl[ (z-s)^{-1}\bigl(\Lambda
_N(z)-s\bigr)^{-1}\biggl(\Lambda
_N(z)-z + \frac{1}{N}g_\sigma(z)^{-1} E_\sigma(z)\biggr) \biggr]\biggr\vert\\
&&\qquad\quad{} + \mathbb{E}\bigl[ \bigl\vert(z-s)^{-1}\bigl\{(z-s)^{-1}-\bigl(\Lambda
_N(z)-s\bigr)^{-1}\bigr\}\bigr\vert\bigr] \frac{1}{N}\vert g_\sigma
(z)^{-1}E_\sigma(z)\vert\\
&&\qquad\leq2\vert\Im m(z)\vert^{-2} \biggl\vert
\Lambda_N(z)-z
+\frac{1}{N}E_\sigma(z)g_\sigma(z)^{-1}\biggr\vert\\
&&\qquad\quad{} + \frac{P_8(\vert
\Im m(z)\vert^{-1})}{N}\vert
\Lambda_N(z) - z\vert(\vert z \vert+K),
\end{eqnarray*}
where we made use of the estimates (\ref{gmoins}), (\ref{binfLambda}),
$ \forall z \in{\mathbb C}\setminus{\mathbb R},$
%
\begin{eqnarray} \label{E}
\forall x \in{\mathbb R} \qquad \biggl\vert\frac{1}{ z-x }
\biggr\vert&\leq&
|\Im m(z)|^{-1}, \nonumber\\[-8pt]\\[-8pt]
\vert E_\sigma(z) \vert&\leq& P_4(\vert\Im m(z)\vert^{-1})\qquad
\mbox{[using (\ref{0.1})]}.\nonumber
\end{eqnarray}
Let us write
\begin{eqnarray*}
& & \biggl\vert\Lambda_n(z)-z
+\frac{1}{N}g_\sigma(z)^{-1}\biggr\vert\\
&&\qquad =  \frac{1}{g_N(z)} \biggl( \sigma^2 g_N^2(z) -z g_N(z) + 1 + \frac
{E_\sigma(z)}{N}
\biggr)\\
&&\qquad\quad{} +\frac{ E_\sigma(z)/N}{g_N(z)g_\sigma(z)} \bigl(g_N(z)-g_\sigma(z)\bigr).
\end{eqnarray*}
We get from Theorem
\ref{theomasteq2}, (\ref{gNmoins}), (\ref{g-gN}), (\ref{E}),
(\ref{gmoins}),
\[
\biggl\vert\Lambda_N(z)-z
+\frac{1}{N}E_\sigma(z)g_{\sigma}(z)^{-1}\biggr\vert\leq(| z|+K)^3\frac
{P_{15}(|\Im m (z)|^{-1})}{N^2}.
\]
Finally, using also (\ref{lambdaN}), we get
for any
$z$ in
$\mathcal{O}'_n$,
\[
\biggl\vert g_{\sigma}(z)-g_N(z) +\frac{1}{N} L_{\sigma}(z) \biggr\vert\leq(|
z|+K)^3\frac{P_{17}(|\Im m (z)|^{-1})}{N^2}.
\]
Now, for $z \notin\mathcal{O}'_n$, such that $\Im m (z) >0$,
\begin{eqnarray*}
1 &\leq& 4\frac{P_6(\vert(\Im m(z))^{-1}\vert)}{N}( \vert z \vert+
\sigma^2 \vert
\Im m(z)\vert^{-1})\vert
\Im m(z)\vert^{-1} \\
&\leq& (\vert z \vert+K)
\frac{P_{8}(\vert\Im m(z)^{-1}\vert)}{N}.
\end{eqnarray*}
We get
\begin{eqnarray*}
& & \biggl\vert g_{\sigma}(z)-g_N(z) +\frac{1}{N} L_{\sigma}(z) \biggr\vert\leq
\vert g_{\sigma}(z)-g_N(z)\vert+ \frac{1}{N}\vert L_\sigma(z)
\vert\\
&&\qquad\leq (\vert z \vert+K) \frac{P_{8}(\vert
\Im m (z)\vert^{-1})}{N} \\
&&\qquad\quad{} \times\biggl[ (\vert z
\vert+K) \frac{P_{9}(\vert\Im m(z)\vert^{-1})}{N} +
\frac{1}{N}P_7(\vert
\Im m(z)\vert^{-1}) (\vert z \vert+K) \biggr]\\
&&\qquad\leq (\vert z \vert+K)^2 \frac{P_{17}(\vert
\Im m(z)\vert^{-1})}{N^2}.
\end{eqnarray*}
Thus,
for any $z$ such that $\Im m (z)>0$,
%
\begin{equation}\label{estimdifL}
\biggl\vert g_{\sigma}(z)-g_N(z)
+\frac{1}{N} L_{\sigma}(z) \biggr\vert\leq(\vert z \vert+K)^3
\frac{P_{17}(\vert\Im m(z)\vert^{-1})}{N^2}.
\end{equation}
Let us denote for a while $g_N=g_N^{A_N}$ and $L_{\sigma
}=L_{\sigma}^{A_N}$.
Note that we get exactly the same estimation (\ref{estimdifL}) dealing
with $-A_N$
instead of $A_N$. Hence since $g_{\sigma}(z)=- g_{\sigma}(-z)$,
$g_N^{-A_N}(z)=-g_N^{A_N}(-z)$ (using the symmetry assumption on $\mu$)
and $L_{\sigma}^{-A_N}(z)=L_{\sigma}^{A_N}(-z)$,
it readily follows that (\ref{estimdifL}) is also valid for
any $z$ such that $\Im m (z)<0$. In conclusion:
\begin{proposition} \label{propestimdifL2}
For any $z \in{\mathbb C}\setminus{\mathbb R}$,
%
\begin{equation} \label{estimdifL2}
\biggl\vert g_{\sigma}(z)-g_N(z)
+\frac{1}{N} L_{\sigma}(z) \biggr\vert\leq(\vert z \vert+K)^3
\frac{P_{17}(\vert\Im m(z)\vert^{-1})}{N^2}.
\end{equation}
\end{proposition}

\subsection{The spectrum of $M_N$}\label{incsp}
The following step now consists of deducing Proposition \ref{epf}
from Proposition \ref{propestimdifL2} (from which we will easily
deduce the appropriate inclusion of the spectrum of $M_N$). Since
this transition is based on the inverse Stieltjes transform, we
start with establishing the fundamental Proposition~\ref{LSt} below
concerning the nature of $L_{\sigma}$. To this aim, it
will be relevant to rewrite $L_{\sigma}$ as
%
\begin{equation}\label{Lsigma2}\qquad\quad
L_{\sigma}(z)= g_{\sigma}(z) ^{-1} \times g_{\sigma}'(z) \times
\Biggl( \sum_{j=1}^J k _j   \frac{\theta
_j}{({1}/({g_{\sigma}(z)}))- \theta_j}+ \frac{\kappa_4}{2}
g_{\sigma}^4(z)\Biggr).
\end{equation}
We recall that $J_{+ \sigma}$ (resp., $J_{- \sigma}$) denotes the
number of $j$'s such that $\theta_j > \sigma$ (resp., $\theta_j <
-\sigma$). As in the \hyperref[intro]{Introduction}, we define
\begin{eqnarray*}{\label{defrhotheta1}}
\rho_{\theta_j}=\theta_j + \frac{\sigma^2}{\theta_j}
\end{eqnarray*}
which is $> 2 \sigma$ (resp., $<-2 \sigma$) when $\theta_j> \sigma$
(resp., $<-\sigma$).

\begin{proposition}\label{LSt}
$L_{\sigma}$ is the Stieltjes transform of a distribution $\Lambda
_{\sigma}$ with
compact support
\[
K_{\sigma}(\theta_1, \ldots,\theta_J):=\{\rho_{\theta_J} ;
\ldots
; \rho_{\theta
_{J-J_{-\sigma}+1}}\} \cup
[-2 \sigma,2 \sigma]\cup\{\rho_{\theta
_{J_{+\sigma}}}; \ldots;
\rho_{\theta_1} \} .
\]
\end{proposition}

The proof relies on the following characterization already used in
\cite{S}.

\begin{theorem}[{\cite{T}}]\label{TS}
\begin{itemize}
\item Let $\Lambda$ be a distribution on ${\mathbb R}$ with compact support.
Define the Stieltjes transform of $\Lambda$, $ l\dvtx{\mathbb C}
\setminus{\mathbb R}\rightarrow{\mathbb C}$ by
\[
l(z) =\Lambda\biggl(\frac{1}{z-x}\biggr).
\]
Then $l$ is analytic in ${\mathbb C}
\setminus{\mathbb R}$
and has an analytic continuation to ${\mathbb C}
\setminus\operatorname{supp}(\Lambda)$. Moreover:
\begin{itemize}[(c$_2$)]
\item[(c$_1$)] $l(z) \rightarrow0$ as $|z| \rightarrow\infty,$
\item[(c$_2$)]
there exist a constant $C>0$, an $n \in{\mathbb N}$ and a compact set $K
\subset{\mathbb R}$ containing $\operatorname{supp}(\Lambda)$
such that for any $z \in{\mathbb C}
\setminus{\mathbb R}$,
\[
|l(z)| \leq C \max \{\operatorname{dist}(z,K)^{-n},1\},
\]
\item[(c$_3$)] for any $\phi\in\mathcal{C}^{\infty}({\mathbb
R},{\mathbb R})$ with
compact support
\[
\Lambda(\phi)=-\frac{1}{\pi} \lim_{y \rightarrow0^+} \Im m \int
_{{\mathbb R}}
\phi(x) l(x+iy) \,dx.
\]
\end{itemize}
\item Conversely, if $K$ is a compact subset of ${\mathbb R}$ and if
$l\dvtx {\mathbb C}
\setminus K \rightarrow{\mathbb C}$ is an analytic function
satisfying (c$_1$) and (c$_2$) above, then $l$ is the Stieltjes
transform of a compactly supported distribution $\Lambda$ on ${\mathbb R}$.
Moreover, $\operatorname{supp}(\Lambda)$ is exactly the set of singular points
of $l$ in $K$.
\end{itemize}
\end{theorem}

The following properties of the Stieltjes transform
$g_{\sigma}$ will be useful for showing that $L_{\sigma}$ fulfills
the previous conditions.

\begin{lemma}{\label{lemme-g-sigma}}
$g _{\sigma}$ is analytic and invertible on $\mathbb C \setminus[-2
\sigma,2
\sigma]$ and its
inverse~$z _{\sigma}$ satisfied
\begin{eqnarray*}\label{z-sigma}
z _{\sigma}(g)= \frac{1}{g} + \sigma^2 g\qquad \forall g \in g
_{\sigma}(\mathbb C \setminus[-2 \sigma,2
\sigma]).
\end{eqnarray*}

\begin{longlist}[(a)]
\item[(a)] The complement of the support of $\mu_{\sigma}$ is
characterized as
follows:
\[
x \in\mathbb R \setminus[-2 \sigma,2
\sigma] \quad\Longleftrightarrow\quad\exists g \in\mathbb R ^* \mbox{ such that
}   \biggl\vert
\frac{1}{g} \biggr\vert> \sigma\mbox{ and }   x=z_{\sigma}(g).
\]
\item[(b)]
Given $x \in\mathbb R \setminus[-2 \sigma,2
\sigma]$ and $\theta\in\mathbb R$ such that $\vert\theta\vert>
\sigma$, one has
\[
\frac{1}{g _{\sigma} (x)}= \theta\quad\Longleftrightarrow\quad x=\theta+
\frac{ \sigma^2}{\theta} := \rho_{\theta}.
\]
\end{longlist}
\end{lemma}

This lemma can be easily proved using, for example, the
explicit expression of~$g_{\sigma}$ [derived from (\ref{g})], namely for all $x \in\mathbb R
\setminus[-2\sigma,2\sigma]$,
\[
g_{\sigma}(x) = \frac{x}{2 \sigma^2} \bigl( 1 - \sqrt{1-{4
\sigma^2}/{x^2}} \bigr).
\]

\begin{pf*}{Proof of Proposition \ref{LSt}}
Using (\ref{Lsigma2}), one readily sees that the set of singular points of
$L_{\sigma}$ is $[-2\sigma, 2 \sigma] \cup\{x \in
\mathbb R \setminus[-2 \sigma,2 \sigma],
\frac{1}{g_{\sigma}(x)} \in\operatorname{Spect}(A_N)\} $.
Hence [using point (b) of Lemma \ref{lemme-g-sigma}] the set of
singular points of $L_{\sigma}$ is exactly $K_{\sigma}(\theta_1,
\ldots
,\theta_J).$

Now, we are going to show that $L_{\sigma}$ satisfies
(c$_1$) and (c$_2$) of Theorem \ref{TS}.
We have obviously that
\[
\vert z- \sigma^2 g_{\sigma}(z)-\theta_j\vert\geq\bigl\vert\vert z
-\theta_j\vert-\vert\sigma^2 g_{\sigma}(z)\vert\bigr\vert.
\]
Now,
let $\alpha> 0$ such that $\alpha> 2 \sigma$ and for any $j=1,
\ldots,
J$, $\alpha- \vert\theta_j \vert>
\frac{\sigma^2}{\alpha- 2 \sigma}$. For any $z \in\mathbb{C}$ such
that $|z|> \alpha$,
\[
\vert z -\theta_j\vert\geq\vert z\vert-\vert\theta_j\vert> \frac
{\sigma^2}{\alpha- 2 \sigma}
\]
and according to (\ref{bgout})
\[
\vert\sigma^2 g_{\sigma}(z) \vert\leq\frac{\sigma^2}{\vert z
\vert-
2 \sigma}\leq\frac{\sigma^2}{\alpha- 2 \sigma}.
\]
Thus we get that for $z \in\mathbb{C}$ such that $|z|>
\alpha$,
\[
\vert z- \sigma^2 g_{\sigma}(z)-\theta_j\vert\geq\vert z\vert
-\vert
\theta_j\vert- \frac{\sigma^2}{\alpha- 2 \sigma}.
\]
Using also (\ref{bgout})--(\ref{bgmoinsout}), we
get readily that for $|z|> \alpha$,
\begin{eqnarray*}
\vert L_{\sigma}(z) \vert&\leq&
\biggl(|z|+ \frac{\sigma^2}{|z| - 2 \sigma}\biggr)
\frac{1}{(|z| - 2 \sigma)^2 }\\
&&{}\times\Biggl(\sum_{j=1}^J \frac{k_j \vert
\theta_j \vert}{ \vert z\vert- \vert\theta_j\vert -
({\sigma^2}/({\alpha- 2 \sigma}))} + \frac{\vert\kappa_4
\vert}{2(|z| - 2 \sigma)^4 } \Biggr).
\end{eqnarray*}
Then, it is clear that $\vert L_{\sigma}(z) \vert\rightarrow0$ when
$|z| \rightarrow+ \infty$ and (c$_1$) is satisfied.

Now we follow the approach of \cite{S} (Lemma 5.5) to prove (c$_2$).
Denote by $\mathcal{E}$ the convex envelope of $K_{\sigma}(\theta_1,
\ldots
,\theta_J)$
and define the interval
\begin{eqnarray*}
K  :\!&= & \{x \in{\mathbb R}; \operatorname{dist}(x, \mathcal{E})\leq1\}
\\
& = & [ \min\{x \in K_{\sigma}(\theta_1, \ldots,\theta_J)\}
-1; \max
\{x \in K_{\sigma}(\theta_1, \ldots,\theta_J)\}+1 ]
\end{eqnarray*}
and
\[
D= \{ z \in{\mathbb C}; 0<\operatorname{dist}(z,K) \leq1\}.
\]
\begin{itemize}
\item Let $z \in D \cap{\mathbb C}\setminus{\mathbb R}$ with $\Re
e(z) \in K$. We
have $\operatorname{dist}(z, K)= |\Im m (z)| \leq1$.
Using the upper bounds
(\ref{bornesupg}), (\ref{gmoins}), (\ref{bgprime}) and
(\ref{0.1}),
we easily deduce that there exists some constant $C_0$ such that for
any $z \in D \cap{\mathbb C}\setminus{\mathbb R}$ with $\Re e(z) \in
K$,
\[
\vert L_{\sigma}(z) \vert\leq C_0 |\Im m (z) |^{-7}= C_0 \operatorname{dist}(z,
K)^{-7}=C_0\max(\operatorname{dist}(z, K)^{-7};1).
\]
\item Let $z \in D \cap{\mathbb C}\setminus{\mathbb R}$ with $\Re
e(z) \notin K$.
Then $\operatorname{dist} (z, K_{\sigma}(\theta_1, \ldots,\theta_J)) \geq1$.
Since $L_{\sigma}$ is bounded on compact subsets of ${\mathbb
C}\setminus
K_{\sigma}(\theta_1, \ldots,\theta_J)$,
we easily deduce that there exists some constant $C_1$ such that for
any $z \in D$ with $\Re e(z) \notin K$,
\[
\vert L_{\sigma}(z) \vert\leq C_1 \leq C_1 \operatorname{dist}(z, K)^{-7}= C_1\max
(\operatorname{dist}(z, K)^{-7};1).
\]
\item Since $\vert L_{\sigma}(z) \vert\rightarrow0$ when $|z|
\rightarrow+ \infty$, $L_\sigma$ is bounded on ${\mathbb
C}\setminus D$.
Thus, there exists some constant $C_2$ such that for any $z \in
{\mathbb C}
\setminus D$,
\[
\vert L_{\sigma}(z) \vert\leq C_2 = C_2\max(\operatorname{dist}(z, K)^{-7};1).
\]
\end{itemize}
Hence (c$_2$) is satisfied with $C= \max(C_0, C_1, C_2)$ and $n=7$ and
Proposition \ref{LSt} follows from Theorem \ref{TS}.
\end{pf*}

We are now in position to deduce the following proposition
from the estimate~(\ref{estimdifL2}).

\begin{proposition}\label{epf}
For any smooth function $\varphi$ with compact support,
%
\begin{equation}\label{esf}
\mathbb{E}[ \mathrm{tr}_N (\varphi(M_N))] = \int\varphi\, d\mu_{sc}+
\frac{1}{N}\Lambda_{\sigma}(\varphi) + O\biggl(\frac{1}{N^2}\biggr).
\end{equation}
Consequently,
for $\varphi$ smooth, constant outside a compact set and such that
$\operatorname{supp}(\varphi) \cap K_{\sigma}(\theta_1, \ldots,\theta_J) =
\varnothing$,
%
\begin{equation}\label{vide} \mathrm{tr}_N (\varphi
(M_N))=O({N^{-{4}/{3}}})\qquad
\mbox{a.s.}
\end{equation}
\end{proposition}

\begin{pf} Using the inverse Stieltjes transform, we
get, respectively, that, for any $\varphi$ in $\mathcal{C}^{\infty
}({\mathbb R}, {\mathbb R})$ with compact support,
\[
\mathbb{E}[ \mathrm{tr}_N (\varphi
(M_N))] - \int
\varphi\, d\mu_{sc}-\frac{\Lambda_{\sigma}(\varphi)}{N}= -\frac
{1}{\pi}
\lim_{y
\rightarrow0^+} \Im m \int_{{\mathbb R}} \varphi(x)
r_N(x+iy)\,dx,
\]
where $r_N= g_\sigma(z)-g_N(z) +\frac{1}{N} L_\sigma(z)$ satisfies,
according to Proposition \ref{propestimdifL2}, for any
$z \in{\mathbb C}\setminus{\mathbb R}$,
\[
\vert r_N(z)\vert\leq\frac{1}{N^2}(\vert z \vert+K)^{\alpha}
P_{k}(\vert\Im m(z)^{-1} \vert),
\]
where $\alpha=3$
and $k=17$.

We refer the reader to the \hyperref[app]{Appendix} of \cite{CD} where it
is proved using the ideas of \cite{HT} that
\[
\limsup_{y \rightarrow0^+} \biggl\vert\int_{{\mathbb R}}
\varphi(x) h(x+iy)\,dx\biggr\vert\leq C< + \infty,
\]
when
$h$ is an analytic function on ${\mathbb C}\setminus{\mathbb R}$
which satisfies
\[
\vert h(z)\vert\leq(\vert z \vert+K)^{\alpha} P_{k}(\vert\Im
m(z)^{-1} \vert).
\]
Dealing with $h(z) =N^2
r_N(z)$, we deduce that
\[
\limsup_{y \rightarrow0^+}
\biggl\vert
\int_{{\mathbb R}} \varphi(x)
r_N(x+iy)\,dx\biggr\vert\leq\frac{C}{N^2}
\]
and then (\ref{esf}).

Following the proof of Lemma 5.6 in
\cite{S}, one can show that $\Lambda_{\sigma}(1)=0$. Then,
the rest of the proof of (\ref{vide}) sticks to the proof of Lemma 6.3
in \cite{HT} (using Lemma \ref{variance}).
\end{pf}

Following \cite{HT} (Theorem 6.4), we set $K=K_{\sigma
}(\theta_1,
\ldots,\theta_J)+(-\frac{\varepsilon}{2},\frac{\varepsilon}{2})$,
$F=\{t \in{\mathbb R}; \operatorname{dist}(t, K_{\sigma}(\theta_1, \ldots,\theta_J))
\geq
\varepsilon\}$ and take $\varphi\in{\mathcal{\mathbb C}}^{\infty
}({\mathbb R},{\mathbb R})$ such
that $0 \leq\varphi\leq1$, $\varphi(t) = 0$ for $t \in K$ and
$\varphi(t) = 1$ for $t \in F$. Then according to
(\ref{vide}), $\mathrm{tr}_N (\varphi(M_N))
=O({N^{-{4}/{3}}}) $ a.s. Since $\varphi\geq1_F$, it
follows that $\mathrm{tr}_N
(1_F(M_N))=O({N^{-{4}/{3}}})$ a.s. and thus the number of
eigenvalues of $M_N$ in $F$ is almost surely an
$O({N^{-{1}/{3}}})$ as $N$ goes to infinity. Since for
each $N$ this number has to be an integer we deduce that the number of
eigenvalues of $M_N$ in $F$ is zero almost surely as $N$ goes to infinity.
The
fundamental inclusion (\ref{inclusion}) follows, namely, for any
$\varepsilon>0$,
almost surely
\[
\operatorname{Spect}(M_N) \subset K_{\sigma}(\theta_1, \ldots,\theta
_J)+(-\varepsilon,\varepsilon),
\]
when $N$ goes to infinity.

Such a method can be carried out in the case of Wigner real
symmetric matrices; then the approximate master equation is the
following [compare with (\ref{masteq})]:
\begin{eqnarray*}
&&\sigma^2 g_N(z)^2 -z g_N(z) +1 +
\frac{1}{N}\frac{\kappa_4}{2}\mathbb{E}\Biggl[\Biggl(\frac{1}{N}\sum_{k=1}^N
G_{kk}(z)^2\Biggr)^2\Biggr] + \frac{\sigma^2}{N}
\mathbb{E}(\mathrm{tr}_N G_N(z)^2) \\
&&\qquad{}+ \mathbb{E}(\mathrm{tr}_N
[A_N G_N(z)])= O\biggl(\frac{1}{N^2}\biggr).
\end{eqnarray*}
Note that the additional term $\frac{\sigma^2}{N}
\mathbb{E}(\mathrm{tr}_N G_N(z)^2)$ already appears in the
non-Deformed GOE case in \cite{S}. One can establish in a similar way
the analogue of \textup{(\ref{g-gN})} and then, following the proof
of Corollary
3.3 in \cite{S}, deduce that
\[
\mathbb{E}(\mathrm{tr}_N G_N(z)^2)= \mathbb{E}
\bigl((z-s)^{-2}
\bigr)+O\biggl(\frac{1}{N}\biggr),
\]
where $s$ is a centered semicircular variable with
variance $\sigma^2$. Hence by similar arguments as in the complex
case, one gets the master equation
\[
\sigma^2 g_N(z)^2 -z g_N(z) +1 + \frac{1}{N}E_{\sigma}(z) =
O\biggl(\frac{1}{N^2}\biggr),
\]
where
\[
E_{\sigma}(z) =\sum_{j=1, j \not= j_0}^J k_j \frac{ \theta_j}{ z -
\sigma^2 g_\sigma(z) - \theta_j} + \frac{\kappa_4}{2} g_\sigma^4(z)
+\mathbb{E}\bigl((z-s)^{-2}\bigr).
\]
It can be proved that $L_\sigma(z):= g_\sigma(z)^{-1} \mathbb
{E}((z-s)^{-2})E_\sigma(z)$ is the Stieltjes transform of a
distribution $\Lambda_\sigma$ with
compact support $K_{\sigma}(\theta_1, \ldots,\theta_J)$, too. The last
arguments hold likewise in the real symmetric case.

Hence we have established:

\begin{theorem}{\label{ThmNooutside}} Let $(M_N)_N$ be any real or
complex Deformed model satisfying \textup{(i)} and
\textup{(ii)} in Section \ref{results}. Let $J_{+ \sigma}$ (resp., $J_{-
\sigma}$) be the number of j's such that $\theta_j > \sigma$ (resp.,
$\theta_j < -\sigma$). Then for any $\varepsilon>0$, almost surely,
there is no eigenvalue of $M_N$ in
%
\begin{eqnarray}\label{reunionoutside}\qquad\quad
&&(-\infty,\rho_{\theta_J} - \epsilon)
\cup(\rho_{\theta_J} + \epsilon,
\rho_{\theta_{J-1}} - \epsilon) \cup\cdots\cup(\rho_{\theta
_{J-J_{-\sigma}+1}}+ \epsilon, -2\sigma- \epsilon)
\nonumber\\[-8pt]\\[-8pt]
&&\qquad{}  \cup  (2\sigma+ \epsilon, \rho_{\theta
_{J_{+\sigma}}}-\epsilon) \cup\cdots\cup(\rho_{\theta_2} +
\epsilon,\rho_{\theta_1}- \epsilon) \cup(\rho_{\theta_1} +
\epsilon,+\infty),\nonumber
\end{eqnarray}
when $N$ is large enough.
\end{theorem}

\begin{remark}
As soon as $\epsilon>0$ is small enough, the union
(\ref{reunionoutside}) is made of nonempty disjoint intervals.
\end{remark}

\subsection{The almost sure convergence result}\label{asc}

As announced in the \hyperref[intro]{Introduction}, Theorem \ref{ThmASCV} is the analogue
of the main statement of \cite{BS3} established for general spiked
population models (\ref{spike}). The previous Theorem \ref
{ThmNooutside} is the main step of the proof since now, we can adapt
the arguments needed for the conclusion of \cite{BS3} viewing the
Deformed Wigner model (\ref{Dw}) as the additive analogue of the spiked
population model (\ref{spike}).

Let us consider one of the positive eigenvalues $\theta_j$ of the
$A_N$'s. We recall that this implies that $\lambda_{k_1+ \cdots
+k_{j-1} +i}(A_N)= \theta_j$ for all $1 \leq i \leq k_j$. We want
to show that if $\theta_j>\sigma$ (i.e., with our notation, if $j
\in\{1, \ldots, J_{+ \sigma}\}$), the corresponding eigenvalues of
$M_N$ almost surely jump above the right endpoint $2\sigma$ of the
semicircle support as
\[
\forall1 \leq i \leq k_j\qquad  \lambda_{k_1+ \cdots+k_{j-1} +i}(M_N)
\longrightarrow\rho_{\theta_j}\qquad \mbox{a.s.},
\]
whereas the rest of the asymptotic spectrum of $M_N$ lies below $2
\sigma$ with
\[
\lambda_{k_1+ \cdots+k_{J_{+ \sigma}} +1}(M_N)
\longrightarrow2\sigma\qquad\mbox{a.s.}
\]
Analogous results hold for
the negative eigenvalues $\theta_j$ [see points (c) and (d) of
Theorem \ref{ThmASCV}]. To describe the phenomenon, one can say
that, when $N$ is large enough, the (first extremal) eigenvalues of
$M_N$ can be viewed as a ``smoothed'' deformation of the (first
extremal) eigenvalues of $A_N$. According to the analysis made in
the previous section [Lemma \ref{lemme-g-sigma}(b)], we already know
that the limits $\rho_{\theta
_j}$ are related to the $\theta_j$'s through the Stieltjes transform
$g_{\sigma}$. More precisely, one has
\[
\mbox{for all $j $ such
that $|\theta_j| > \sigma$,}\qquad
\frac{1}{g_{\sigma}(\rho_{\theta_j})}= \theta_j .
\]
Our main
purpose now is to establish the asymptotic link between the spectra
of the matrices $M_N=X_N+A_N$ and $A_N$.

Intuitively, this link seems rather natural when $\sigma$ is close
to zero. Indeed, when~$N$ goes to infinity, since the spectrum of
$X_N$ is concentrated in $[-2\sigma,2\sigma]$
[recall (\ref{extremal})], the spectrum of $M_N$ should be close to
the one of $A_N$ as soon as $\sigma$ will be close to zero (in other
words, the spectrum of $M_N$ is, viewed as a deformation of the one
of $A_N$, continuous in $\sigma$ in a neighborhood of zero). Thus
given an interval $[a,b] \subset{} {^{c}\!K}_{\sigma}(\theta_1, \ldots,
\theta_J)$, the result of Theorem \ref{ThmNooutside} saying that
$[a,b]$ does not contain eigenvalues of $M_N$ should be improved: it
should correspond to $[a,b]$ some
interval $I$ close to $[a,b]$, lying outside the spectrum of $A_N$ and
such that the number of eigenvalues of
$M_N$ in one side of $[a,b]$ is equal to the one of $A_N$ in the
corresponding side of $I$. Following \cite{BS2}, we will say that
there is exact separation of eigenvalues of the matrices $A_N$
and $M_N$.

In the following section, we justify that the exact
separation phenomenon occurs regardless of the size of $\sigma$. The
proof of Theorem \ref{ThmASCV} will then follow from some suitable
choices of $[a,b]$ (see Section
\ref{sectionProoffin}).

\subsubsection{Exact separation of eigenvalues}{\label{sectionExact}}

According to the previous discussion, we need to refine the analysis
made on $g_{\sigma}$ in order to identify and understand the link
between intervals in $ ^c K_{\sigma}(\theta_1, \ldots,
\theta_J)$ and the complement of the spectrum
of the $A_N$'s. We also need to understand the dependence on $\sigma$.
This is the aim of the following important Lemma \ref{lemme-g-sigma2}.

As before,
we denote (recall Lemma \ref{lemme-g-sigma}) by $z_{\sigma}$ the
inverse function of $g_{\sigma}$ which is given by
\[
z_{\sigma}(g)=\frac{1}{g}+ \sigma^2 g.
\]
Using Lemma \ref{lemme-g-sigma}, one readily sees that the set
$ ^cK_{\sigma}(\theta_1, \ldots, \theta_J)$ can be characterized
as follows:
%
\begin{equation} {\label{Complem-SuppL}}
x \in {} {^{c}\!K}_{\sigma}(\theta_1, \ldots,
\theta_J)\quad\Longleftrightarrow\quad
\exists g \in
\mathcal G _{\sigma} \mbox{ such that }x= z_{\sigma}(g),
\end{equation}
where
%
\begin{equation}\label{Gsigma}
\mathcal G _{\sigma}:= \biggl\{ g
\in
\mathbb R ^*,   \biggl|\frac{1}{g}\biggr| >
\sigma\mbox{ and } \frac{1}{g} \notin\operatorname{Spect}(A_N) \biggr\}.
\end{equation}
Obviously, one has $g=g_{\sigma}(x)$ if $x \in {} {^{c}\!K}_{\sigma}(\theta
_1, \ldots, \theta_J)$.

\begin{lemma} {\label{lemme-g-sigma2}}
Let $[a,b]$ be a compact set contained in $ ^c K_{\sigma}(\theta_1,
\ldots, \theta_J)$. Then:
\begin{longlist}
\item ${[\frac{1}{g_{\sigma}(a)}, \frac
{1}{g_{\sigma
}(b)}]} \subset(\operatorname{Spect}(A_N))^c$.
\item  For all $0 < \hat{\sigma} < \sigma$, the interval
$[z_{\hat
{\sigma}}(g_{\sigma}(a)),z_{\hat{\sigma}}(g_{\sigma}(b))]$ is contained
in $ ^c K_{\hat{\sigma}}(\theta_1, \ldots, \theta_J)$ and $z_{\hat
{\sigma}}(g_{\sigma}(b))-z_{\hat{\sigma}}(g_{\sigma}(a)) \geq b-a$.
\end{longlist}
\end{lemma}

\begin{pf}
The function $1/{g _{\sigma}}$ being increasing, (i) readily follows
from (\ref{Complem-SuppL}).

Noticing that $\mathcal G _{{\sigma}} \subset\mathcal G
_{\hat{\sigma}}$ for all $\hat{\sigma} < \sigma$ implies (recall
also that $g_{\sigma}$ decreases on $[a ,b ]$) that
$[g_{\sigma}(b),g_{\sigma}(a)] \subset\mathcal G _{\hat{\sigma}}$.
Relation (\ref{Complem-SuppL}) combined with the fact that the
function $z_{\hat{\sigma}}$ is decreasing on
$[g_{\sigma}(b),g_{\sigma}(a)]$ leads to
\[
[z_{\hat{\sigma}}(g_{\sigma}(a))
,z_{\hat{\sigma}}(g_{\sigma}(b)) ] \subset {} {^{c}\!
K}_{\hat{\sigma}}(\theta_1, \ldots, \theta_J)
\]
and the first part of (ii) is stated.
Now, we have
\begin{eqnarray*}
l_\sigma(\hat{\sigma}):\!&=& z_{\hat{\sigma}}(g_{\sigma
}(b))-z_{\hat
{\sigma}}(g_{\sigma}(a))\\
&=& \frac{1}{g_{\sigma}(b)}-\frac{1}{g_{\sigma}(a)}+ \hat{\sigma
}^2\bigl(g_{\sigma}(b)-g_{\sigma}(a)\bigr).
\end{eqnarray*}
Since $g_{\sigma}$ decreases on $[a ,b ]$, we have $g_{\sigma
}(b)-g_{\sigma}(a)\leq0$ and thus $
l_\sigma$ is decreasing on ${\mathbb R}^+$. Then the last point of (ii)
follows since $l_\sigma(\sigma)=b-a$.
\end{pf}

The exact separation result can now be stated. Let $[a,b]$ be an
interval contained in $ ^c K_{\sigma}(\theta_1, \ldots,
\theta_J)$. By Theorem \ref{ThmNooutside}, $[a,b]$ is outside the
spectrum of $M_N$. Moreover, from Lemma \ref{lemme-g-sigma2}(i),
there corresponds an interval $I=[a',b']$ outside the spectrum of $A_N$,
that is, there is $i_N \in\{0, \ldots, N \}$ such that
%
\begin{equation}{\label{sep1}}
\lambda_{i_N+1}(A_N) < \frac{1}{g_{\sigma}(a)}:=a'  \quad\mbox{and}\quad
\lambda_{i_N}(A_N) > \frac{1}{g_{\sigma}(b)}:=b'.
\end{equation}
$a$ and $a'$ (resp., $b$ and $b'$) are linked as follows:
\begin{eqnarray*}{\label{lien-a-a'}}
a=\rho_{a'}:= a' + \frac{\sigma^{2}}{a'}  \qquad(\mbox{resp., }
b=\rho_{b'}).
\end{eqnarray*}
We claim that $[a,b]$ splits the eigenvalues of $M_N$ exactly as $I$
splits the spectrum of $A_N$. In other words:

\begin{theorem}{\label{Thmexact}} With $i_N$ satisfying $(\ref
{sep1})$, one has
%
\begin{equation}\label{sep2}
\mathbb P[\lambda_{i_N+1}(M_N) < a   \mbox{ and }
\lambda_{i_N}(M_N) >
b,   \mbox{ for all large $N$}]=1.\\
\end{equation}
\end{theorem}

This result is the analogue of the main statement of \cite{BS2} (cf.
Theorem 1.2 of~\cite{BS2}) established in the spiked population
setting (and in fact for quite general sample covariance matrices).
Its proof is quite technical and is inspired by the work \cite{BS2}.
It mainly relies on results on eigenvalues of the rescaled Wigner
matrix $X_N$ combined with the following classical result (due to
Weyl).

\begin{lemma}[(cf. Theorem 4.3.7 of \protect\cite{HJ})] \label{Weyl}
Let B and C be two $N \times N$ Hermitian matrices. For any pair of integers
$j,k$ such that $1 \leq j,k\leq N$ and $j+k \leq N+1$, we have
\[
\lambda_{j+k-1} (B+C) \leq\lambda_{j}(B) + \lambda_{k}(C).
\]
For any pair of integers
$j,k$ such that $1 \leq j,k\leq N$ and $j+k \geq N+1$, we have
\[
\lambda_j(B) + \lambda_k(C) \leq\lambda_{j+k-N} (B+C).
\]
\end{lemma}

\begin{remark} Note that this lemma is the additive analogue of Lemma 1.1
of~\cite{BS2} needed for the investigations of the spiked population model.
\end{remark}

In particular, Lemma \ref{Weyl} gives that $\lambda_{i_N+1} (M_N)
\leq\lambda_{i_N+1}(A_N) + \lambda_{1}(X_N)$ and $\lambda_{i_N}
(M_N) \geq\lambda_{i_N}(A_N) + \lambda_{N}(X_N).$ Besides, as both
$\lambda_{1}(X_N)$ and $-\lambda_{N}(X_N)$ tend toward $2 \sigma$
as $N \to\infty$ [this is ({\ref{extremal}})], the statement of
Theorem \ref{Thmexact} can be quite easily derived when $\sigma$ is
close enough to zero. To handle the general case, the key idea is
that one can reduce to the previous situation by introducing some
parameters. More precisely, given $L>0$ and $k \geq0$, we will introduce
the Wigner matrix
\[
W_N^{k,L} = W_N / \sqrt{1+{k}/{L}}
\]
and let
\[
M_N^{k,L} =A_N+W_N^{k,L}/ {\sqrt{N}}
\]
be the Deformed Wigner matrix of parameter
\[
\sigma_{k,L}= \sigma/
{\sqrt{1+{k}/{L}}}.
\]
The proof will be organized as
follows. On the one hand, as $\sigma_{k,L} \to0$ when $k \to
\infty$ (for any fixed $L>0$), we will readily prove that exact
separation occurs for the matrices $A_N$ and
$M_N^{K,L}$ as soon as $K$ is large enough. On the other hand, we will
show that exact separation also occurs for the eigenvalues of
$M_N=M_N^{0,L}$ and $ M_N^{K,L}$ choosing $L$ large enough. This latter
point will be established by induction on $k$; the underlying idea is
that when the parameter $L$ is large, the matrices $M_N^{k,L}$ and
$M_N^{k+1,L}$ are close to each other and hence split their spectrum
in a similar way.

\begin{pf*}{Proof of Theorem \protect\ref{Thmexact}} With our choice of
$[a,b]$ and the very definition of the spectrum of the $A_N$'s, one
can consider $\epsilon' >0$ small enough such that, for all large
$N$,
\[
\lambda_{i_N+1}(A_N) < \frac{1}{g_{\sigma}(a)}-\epsilon'
\quad\mbox{and}\quad \lambda_{i_N}(A_N) > \frac{1}{g_{\sigma}(b)}+\epsilon'.
\]
Given $L>0$ and $k \geq0$ (their size will be determined later), we
define
\[
a_{k,L}= z_{\sigma_{k,L}}(g_{\sigma}(a))   \quad\mbox{and}\quad   b_{k,L}=
z_{\sigma_{k,L}}(g_{\sigma}(b)),
\]
where we recall that $z_{\sigma_{k,L}}(g)={1}/{g}+ \sigma_{k,L}^2 g.$
Note that for all $L >0$, one has $a_{0,L}=a$ and $b_{0,L}=b$.

We first choose the size of $L$ as follows. We take $L_0$ large enough
such that for all $L \geq L_0$,
%
\begin{equation}{\label{condL}}
\max\bigl( (\sigma^2/L) \bigl( |g_{\sigma}(a)|+ |g_{\sigma}(b)|\bigr) ;
{3\sigma}/{L} \bigr) < ({b-a})/{4}.
\end{equation}
From the very definition of the $a_{k,L}$'s and $b_{k,L}$'s, one can
easily see that $b_{k,L}-a_{k,L} \geq b-a$ [using the last point of
(ii) in Lemma \ref{lemme-g-sigma2}] and that this choice of $L_0$
ensures that, for all $L\geq L_0$ and for all $k \geq0$,
%
\begin{equation}{\label{condL2}}\qquad
|a_{k+1,L} -a_{k,L}| < ({b-a})/{4}\quad \mbox{and}\quad  |b_{k+1,L}
-b_{k,L}| < ({b-a})/{4}.
\end{equation}
Now, we fix $L$ such that $L \geq L_0$ and we write $a_k=a_{k,L}$,
$b_k=b_{k,L}$ and $\sigma_{k} =\sigma_{k,L}$.

We first show that there exists $K$ large enough such that, for all $k
\geq K$, there is exact separation of the eigenvalues of the matrices
$A_N$ and $M_N^{k,L}$, that is,
%
\begin{equation}{\label{casekgen}} \mathbb P{[}\lambda
_{i_N+1}(M_N^{k,L})<a_{k} \mbox{ and } \lambda_{i_N}(M_N^{k,L})>b_{k}
\mbox{ for all
large $N$}{]}=1.
\end{equation}
Lemma \ref{Weyl} first gives that
\[
\lambda_{i_N+1}(M_N^{k,L}) \leq a_k - \epsilon' - \sigma_k ^2
g_{\sigma} (a) +
\frac{1}{\sqrt{1+{k}/{L}}} \lambda_{1} (X_N) \qquad \mbox{if }i_N<N
\]
and
\[
\lambda_{i_N}(M_N^{k,L}) \geq b_k + \epsilon' - \sigma_k ^2
g_{\sigma} (b) +
\frac{1}{\sqrt{1+{k}/{L}}} \lambda_{N} (X_N) \qquad\mbox{if } i_N>0.
\]
Furthermore, according to ({\ref{extremal}}), the two first
extremal eigenvalues
of $X_N$ are such that almost surely and for all $N $
large enough,
\[
0 <\max( -\lambda_N(X_N),\lambda_1(X_N)) < 3 \sigma.
\]
Thus for all $k$, almost surely, at least for $N $ large enough
($N$ does not depend on~$k$),
\[
0 <\frac{1}{\sqrt{1+{k}/{L}}} \times\max(- \lambda_N(X_N),\lambda
_1(X_N)) < 3 \sigma_k.
\]
As $\sigma_{k} \to0$
when $k\to+\infty$, there is $K$ large enough such that for all $k
\geq K$,
\[
\max\bigl( |3 \sigma_k - \sigma_k ^2 g_{\sigma} (a)|, |3 \sigma_k +
\sigma_k
^2 g_{\sigma} (b)|\bigr) < \epsilon'
\]
and then, almost surely, for all $N$ large enough
%
\begin{equation}\label{inf}\lambda_{i_N+1}(M_N^{k,L})<a_{k}
\qquad\mbox{if  } i_N < N
\end{equation}
and
%
\begin{equation}\label{sup}
\lambda_{i_N}(M_N^{k,L})>b_{k}\qquad \mbox{if  } i_N >0.
\end{equation}
Since $\lambda_{N+1}(M_N^{k,L})=-\lambda
_{0}(M_N^{k,L})=-\infty$, (\ref{inf}) [resp., (\ref{sup})] is obviously
satisfied if $i_N=N $ (resp., $i_N=0$). Thus, we have established that
for any $i_N \in\{0,\ldots,N\}$ satisfying (\ref{sep1}), (\ref
{casekgen}) holds for all $k \geq K$. In particular,
%
\begin{equation}{\label{caseK}}\quad
\mathbb P{[}\lambda_{i_N+1}(M_N^{K,L})<a_{K} \mbox{ and }
\lambda
_{i_N}(M_N^{K,L})>b_{K} \mbox{ for all
large $N$}{]}=1.
\end{equation}

Now, we shall show that with probability $1$: for $N$ large, $[a_{K} ,
b_{K}]$ and $[a,b]$ split the eigenvalues of, respectively, $M_N^{K,L}$
and $M_N$ having equal amount of eigenvalues to the left sides of the
intervals. To this aim, we will proceed by induction on~$k$ and
establish that, for all $ k \geq0$, $[a_{k} , b_{k}]$ and $[a,b]$
split the eigenvalues of~$M_N^{k,L}$ and $M_N$ (recall that
$M_N=M_N^{0,L}$) in exactly the same way. To begin, let us consider for
all $k \geq0$, the set
\[
\mathrm{E} _{k} = \{\mbox{no eigenvalues of $M_N^{k,L}$ in
$[a_{k},b_{k}]$, for all large $N$}
\}.
\]
By Lemma \ref{lemme-g-sigma2}(ii) and Theorem \ref{ThmNooutside}, we
know that
$\mathbb P (\mathrm{E}_{k})=1$ for all $k$. In particular, from the
fact that
$\mathbb P (\mathrm{E}_{0})=1$, one has for all
$\omega\in\mathrm{E}_{0}$ and for all large $N$,
%
\begin{equation}{\label{case0}}\qquad \exists j_N(\omega) \in\{ 0, \ldots
, N
\} \mbox{  such that  }
\lambda_{j_N(\omega)+1}(M_N) < a   \mbox{ and }   \lambda
_{j_N(\omega
)}(M_N) >
b.
\end{equation}
Extending the random variable $j_N$ by setting, for instance,
$j_N:=-1$ on $ ^c\mathrm{E}_{0}$,
we want to show that for all $k$,
%
\begin{equation} {\label{casek}}
\mathbb P[\lambda_{j_N+1}(M_N^{k,L}) < a_{k}   \mbox{ and }
\lambda
_{j_N}(M_N^{k}) >
b_{k},   \mbox{for all large $N$}]=1.
\end{equation}
This can be done by induction calling, one more time, on Lemma \ref
{Weyl}. By (\ref{case0}), this is true for $k=0$. Now, let us assume
that (\ref{casek}) holds true. One has
\[
M_N^{k+1,L}= M_N^{k,L} + \biggl(\frac{1}{\sqrt{1+(k+1)/{L}}} - \frac
{1}{\sqrt
{1+{k}/{L}}}\biggr) X_N
\]
so, by Lemma \ref{Weyl},
\[
\lambda_{j _N+1}(M_N^{k+1,L}) \leq\lambda_{j _N+1}(M_N^{k,L})+
(-\lambda_{N} (X_N))/L.
\]
But, for $N$ large enough, $0< -\lambda_{N} (X_N) \leq
3 \sigma$ a.s., so by the condition (\ref{condL}) on~$L$,
\[
\lambda_{j_N+1}(M_N^{k+1,L}) < a_k + ({b-a})/{4} :=
\hat{a}_k.
\]
Similarly, one can show that
\[
\mbox{a.s.}\qquad \lambda_{j_N}(M_N^{k+1,L}) > b_k - ({b-a})/{4} :=
\hat{b}_k.
\]
By (\ref{condL2}), one readily observes that $\hat{a}_{k}-a _{k+1} =a
_{k}-a _{k+1} + ({b-a})/{4} >0$ and similarly that $\hat{b}_{k}-b
_{k+1}<0$. This implies that
\[
[\hat{a}_k, \hat{b}_k] \subset[{a}_{k+1}, {b}_{k+1}].
\]
As $\mathbb P (\mathrm{E}_{k+1})=1$, we deduce that with probability $1$,
\[
\lambda_{j_N+1}(M_N^{k+1,L}) < {a}_{k+1} \quad  \mbox{and}\quad   \lambda
_{j_N}(M_N^{k+1,L}) > {b}_{k+1}\qquad  \mbox{for all $N$ large}.
\]
As a consequence, ({\ref{casek}}) holds for all $k \geq0$ and in
particular for $k=K$. Comparing this with (\ref{caseK}), we deduce that
$j_N=i_N$ a.s. and
\[
\mathbb P{[}\lambda_{i_N+1}(M_N)<a \mbox{ and } \lambda
_{i_N}(M_N)>b \mbox{ for all
large $N$}{]}=1.
\]
This ends the proof of Theorem \ref{Thmexact}.
\end{pf*}

Now, we are in position to prove the main Theorem \ref{ThmASCV}.

\subsubsection{\texorpdfstring{Proof of Theorem \protect\ref{ThmASCV}}{Proof of Theorem 2.1}}
\label{sectionProoffin}
Our reasoning is close to the last Section 4 of~\cite{BS3}. It is
enough to establish parts (a) and (b)
since the assertions (c) and (d) can then be deduced by taking
$-M_N$ instead of $M_N$.

The proof of (a) is mainly based on successive applications of Theorem
\ref{Thmexact}. Fix an integer $1 \leq j \leq J_{+\sigma}$, and let us
consider for $\epsilon>0$, the interval $[a,b]=[\rho_{\theta
_j}+\epsilon,\rho_{\theta_{j-1}}-\epsilon]$ which is included in the
union (\ref{reunionoutside}) (at least for $\epsilon$ small enough). We
define $K_{j(-1)}=k_1 +
\cdots+k_{j(-1)}$. We also take $\theta_0:=+\infty$ and recall the
conventions that $\lambda_0(M_N)=\lambda_0(A_N)=+\infty$ and
$K_{0}=0$. Since
$1 / g_{\sigma}(\rho_{\theta_{k}})=\theta_{k}$ for $k= j-1$ and $j$
and since the function $1/g_{\sigma}$ is continuous and increasing on
$[a,b]$, the
compact interval $[a,b]$ satisfies (\ref{sep1}) with $i_N=K_{j-1}$.
Hence by Theorem \ref{Thmexact}, one has
\[
\mathbb P[\lambda_{K_{j-1}}(M_N) \geq\rho_{\theta_{j-1}}-\epsilon
\mbox{ and }
\lambda_{K_{j-1}+1}(M_N)\leq\rho_{\theta_{j}}+\epsilon
\mbox{, for $N$ large}]=1.
\]
Similar arguments imply that for all $j \in\{1, \ldots, J_{+\sigma}-1
\}$,
\[
\mathbb P[\lambda_{K_{j}}(M_N)\geq\rho_{\theta_{j}}-\epsilon
\mbox
{ and }
\lambda_{K_{j}+1}(M_N) \leq\rho_{\theta_{j+1}}+\epsilon\mbox{,
for $N$
large}]=1.
\]
As a result, we deduce that for all $1
\leq j \leq J_{+\sigma}-1$,
%
\begin{eqnarray}\label{eqfin}\qquad
&&\mathbb P[\rho_{\theta_{j}} -\epsilon\leq\lambda_{K_j}(M_N) \leq
\cdots\leq\lambda_{K_{j-1}+1}(M_N)\leq\rho_{\theta_{j}}+\epsilon
\nonumber\\[-8pt]\\[-8pt]
&&\hspace*{191pt} \mbox{ for $N$ large}]=1.\nonumber
\end{eqnarray}
So, letting $ \epsilon$ go to zero, we obtain (a) for each
integer $j$ of $\{1, \ldots, J_{+\sigma}-1\}$.

Let us now quickly consider the case where $j=J_{+\sigma}$. Note first
that, from the
preceding discussion, we still have (for $\epsilon$ small enough)
\[
\mathbb P[\lambda_{K_{J_{+\sigma}-1}+1}(M_N)\leq\rho_{\theta
_{J_{+\sigma}}}+\epsilon\mbox{, for $N$
large}]=1.
\]
Then, using the fact that $1/g_{\sigma}$ increases
continuously on
$]2\sigma,+ \infty[$ with $1/\break g_{\sigma}(]2\sigma,+ \infty
[)=\,]\sigma,+
\infty[$, one can show that once $\epsilon>0$ is small enough, the
compact set $[a,b]=[2\sigma+\epsilon, \rho_{\theta_{J_{+\sigma}}}
-\epsilon]$ satisfies the
assumptions of Theorem \ref{Thmexact}
with $i_N=K_{J_{+\sigma}}$. This leads to
\[
\mathbb P[\lambda_{K_{J_{+ \sigma}}}(M_N) \geq\rho_{\theta_{J_{+
\sigma}}} -\epsilon  \mbox{
and }   \lambda_{K_{J_{+ \sigma}}+1}(M_N) \leq2 \sigma+\epsilon
\mbox{, for $N$ large}]=1.
\]
Letting $\epsilon\to0$, we deduce that (\ref{eqfin}) holds for
$j=J_{+ \sigma}$ and the assertion (a) is
established. For point (b), the preceding analysis gives that
$\limsup
_N \lambda_{K_{J_{+ \sigma}}+1}(M_N) \leq2\sigma  \mbox{ a.s.}$ and it
remains to prove that
\[
\liminf_N \lambda_{K_{J_{+ \sigma}}+1}(M_N) \geq2\sigma \qquad \mbox{a.s.}
\]
This inequality follows from the fact that the spectral measure
of $M_N$ converges a.s. toward the semicircle law $\mu
_{sc}$ which is compactly supported in $[-2 \sigma, 2 \sigma]$. This
completes the proof of Theorem \ref{ThmASCV}.

\section{Fluctuations}\label{sec5}
The (complex or real) Wigner matricial models under consideration
are the same as previously [i.e., defined by (i) in Section
\ref{results}]
but now we assume that the perturbation $A_N$ is diagonal: $A_N=
\operatorname{diag}(\theta, 0, \ldots,0)$ with unique nonnull eigenvalue $\theta>
\sigma$.
According to the previous section, the a.s. convergence of
$\lambda_1(M_N)$ toward $\rho_\theta$ is universal in the sense that
it does not depend on
$\mu$.

In the first part of this section, we will show that the
fluctuations of $\lambda_1(M_N)$ around this universal limit are not
universal any more. Indeed, we are going to prove that
$\sqrt{N}(1-{\sigma^2}/{\theta^2})^{-1}(\lambda_1(M_N) - \rho
_{\theta
})$ converges in distribution toward the convolution of $\mu$
and a Gaussian distribution. Hence, the limiting distribution clearly
varies with $\mu$ and in particular cannot be Gaussian unless $\mu$ is
Gaussian.

In the second part of this section, we will sharpen the
analysis of the particular Deformed GOE model and explain how this
gives Theorem \ref{ThmFluctuationsUniversalite}.

\subsection{Basic tools}
We start with the following results which will be of basic use later
on. Note that in the following, a complex random variable $x$ will be
said to be \textit{standardized} if $\mathbb{E}(x)=0$ and $\mathbb{E}(\vert x
\vert^2)=1$.

\begin{theorem}[(Lemma 2.7 \protect\cite{BS1})]\label{BaiSilver98}
Let $B=(b_{ij})$ be an $N \times N$ Hermitian matrix and
$Y_N$ be a vector of size $N$ which contains i.i.d.
standardized entries with bounded fourth
moment. Then there is a constant $K>0$ such that
\[
\mathbb E\vert Y_N^* B Y_N
- \operatorname{Tr}
B\vert^2 \leq K \operatorname{Tr}(BB^*).
\]
\end{theorem}

\begin{theorem}[{(cf. \protect\cite{BY2} or \hyperref[app]{Appendix} by J. Baik
and J. Silverstein)}]\label{Silver}
Let $B=(b_{ij})$ be a $N \times N$ random Hermitian matrix and
$Y_N=(y_1, \ldots, y_N)$ be an independent vector of size $N$ which
contains i.i.d.
standardized entries with bounded fourth
moment and such that $\mathbb E (y_1^2)=0$ if $y_1$ is complex. Assume
that:
\begin{longlist}
\item there exists a constant $a>0$ (not depending on $N$) such
that $\Vert
B\Vert \leq a$,
\item $\frac{1}{N} \operatorname{Tr} B^2$ converges in probability to a
number $a_2$,
\item $\frac{1}{N}\sum_{i=1}^N b_{ii}^2$ converges in probability
to a number $a_1^2$.
\end{longlist}
Then the random variable $({1}/{\sqrt N}) ( Y_N ^* B Y_N -
\operatorname{Tr} B ) $ converges in distribution to a Gaussian
variable with mean zero and variance
\[
( \mathbb E |y_1|^4 -1-t/2) a_1 ^2 + (t/2) a_2,
\]
where $t=4$ when $y_1$ is real and is $2$ when $y_1$ is complex.
\end{theorem}

\begin{pf} This result is in fact a particular case of a more general
result of \cite{BY2} (Theorems 7.1 and 7.2) which follows from the method
of moments. We give an alternative elegant proof by J. Baik and J. Silverstein
in the \hyperref[app]{Appendix} of the present paper.
\end{pf}

\begin{theorem}[(Theorem 1.1 in \protect\cite{BY1})]\label{BaiYao}
Let $f$ be an analytic function on an open set of the complex plane
including $[-2 \sigma, 2\sigma]$. If the entries of a general Wigner
matrix $W_N=((W_N)_{ij})_{1 \leq i \leq j \leq N}$ satisfy the
conditions:
\begin{itemize}
\item for $i \neq j$, ${\mathbb E}(\vert(W_N)_{ij} \vert^4) = const$,
\item for any $\eta> 0$, $\lim_{ N \rightarrow+ \infty} \frac
{1}{\eta
^4 N^2} \sum_{i,j}
{\mathbb E}[ \vert(W_N)_{ij} \vert^4 1_{\{\vert
(W_N)_{ij}\vert\geq\eta\sqrt{N}\}}]=0,$
\end{itemize}
then the random variable $ N ( \mathrm{tr}_N(f(\frac{1}{\sqrt N}
W_N))-\int f \,d \mu_{sc} )$ converges in distribution toward
a Gaussian variable.
\end{theorem}

In our setting, $\mu$ satisfies a Poincar\'e inequality and thus, as
already noticed in Section \ref{results}, $\mu$ satisfies
$\int\vert x \vert^q \,d \mu(x) < + \infty$ for any $q$ in ${\mathbb N}$.
Hence, the general Wigner matrices we consider obviously satisfy the
conditions of Theorem \ref{BaiYao}.
Nevertheless, in the following study of fluctuations, we do not use the
Poincar\'e inequality; thus one can expect
that Theorem \ref{ThmFluctuations} is still valid under assumptions on
the only four first moments of $\mu$ provided
one can prove the a.s. convergence of $\lambda_1(M_N)$ toward $\rho
_\theta$
under these weaker assumptions.

\subsection{\texorpdfstring{Proof of Theorem
\protect\ref{ThmFluctuations}}{Proof of Theorem 2.2}}
\label{sectionThmfluctu}
The approach is the same for the complex and real settings and is
close to the one of \cite{Pa} and the ideas of \cite{BBPbis}. Let
$\widehat{M}_{N-1}$ be the $N-1 \times N-1$ matrix obtained from
$M_N$ removing the first row and the first column. Thus,
$\sqrt{N/(N-1)} \widehat{M}_{N-1}$ is a non-Deformed
Wigner matrix associated with the measure $\mu$.
We denote by $\lambda_1(\widehat{M}_{N-1})$ [resp., $\lambda
_{N-1}(\widehat{M}_{N-1})$] the largest (resp., lowest) eigenvalue of $
\widehat{M}_{N-1}$.

Let $0< \delta< {(\rho_{\theta}-2\sigma)}/{4}$. Let us define the
event
\[
\Omega_N= \{ \lambda_1(\widehat{M}_{N-1}) \leq2 \sigma
+ \delta; \lambda_{N-1}(\widehat{M}_{N-1}) \geq-2 \sigma- \delta;
\lambda_1({M}_{N}) \geq\rho_{\theta}- \delta\}.
\]
According to ({\ref{extremal}}) and Theorem \ref{ThmASCV},
$\lim_{N \rightarrow+ \infty} \mathbb P(\Omega_N)= 1.$
Thus, it is sufficient to restrict ourselves to the event $\Omega_N$
in order to study the convergence in distribution of $\sqrt
{N}(1-{\sigma
^2}/{\theta^2})^{-1}(\lambda_1(M_N) - \rho_{\theta})$.

Let $V=\, {^t}\! (
v_1, \ldots, v_N
)$ be an eigenvector corresponding to $\lambda_1(M_N)$. Define the
following vectors in ${{\mathbb C}}^{N-1}$:
\[
\widehat{V}= \, {^t}\! (
v_2 ,\ldots, v_N
)
\]
and
\[
\check{M}_{\bolds\cdot1}=  \, {^t}\! (
(M_N)_{21}, \ldots, (M_N)_{N1})= \frac{1}{\sqrt{N}}  {^t}\! (
(W_N)_{21}, \ldots, (W_N)_{N1}).
\]
Then,
\[
M_N V = \lambda_1(M_N)V \quad\Longleftrightarrow\quad
\cases{\displaystyle
\theta v_1 + \frac{(W_N)_{11}}{\sqrt{N}}v_1+ \check{M}_{\bolds\cdot1}^*
\widehat{V} = \lambda_1(M_N)v_1, \cr
\check{M}_{\bolds\cdot1}v_1 + \widehat{M}_{N-1} \widehat{V}= \lambda_1(M_N)
\widehat{V}.}
\]
On $\Omega_N$, $\lambda_1(M_N)$ is not an eigenvalue of
$\widehat
{M}_{N-1}$ and one can write the eigen-equations
using the resolvent $\widehat{G}(\lambda_1(M_N)):= (\lambda
_1(M_N)I_{N-1} -
\widehat{M}_{N-1})^{-1}$ as follows:
%
\begin{eqnarray}
\widehat{V}&=& v_1 \widehat{G}(\lambda
_1(M_N))\check
{M}_{\bolds\cdot1},\\
\lambda_1(M_N)v_1&=&\theta v_1 + \frac{(W_N)_{11}}{\sqrt{N}}v_1+ v_1
\check
{M}_{\bolds\cdot1}^*\widehat{G}(\lambda_1(M_N))\check{M}_{\bolds\cdot1}.
\label{un}
\end{eqnarray}
Since $v_1$ is obviously nonequal to zero, one gets from (\ref{un})
%
\begin{equation}\lambda_1(M_N)=\theta+ \frac{(W_N)_{11}}{\sqrt{N}}+
\check
{M}_{\bolds\cdot1}^*\widehat{G}(\lambda_1(M_N))\check{M}_{\bolds\cdot1}.\label{deux}
\end{equation}
Moreover, on $\Omega_N$, $\rho_{\theta}$ is not an eigenvalue of
$\widehat{M}_{N-1}$ (recall that $\rho_{\theta} > 2 \sigma$) and the
resolvent $\widehat{G}(\rho_{\theta}):= (\rho_{\theta} I_{N-1} -
\widehat{M}_{N-1})^{-1}$ is well defined, too.
Thus, (\ref{deux}) is equivalent to
\[
\lambda_1(M_N)-\rho_{\theta} = \frac
{(W_N)_{11}}{\sqrt{N}}+
\check{M}_{\bolds\cdot1}^*\widehat{G}(\rho_{\theta})\check{M}_{\bolds\cdot1}-
\frac{\sigma^2}{\theta} +\check{M}_{\bolds\cdot1}^*[\widehat
{G}(\lambda_1(M_N))-
\widehat{G}(\rho_{\theta})]\check{M}_{\bolds\cdot1}.
\]
Using
$ \widehat{G}(\lambda_1(M_N))-\widehat{G}(\rho_{\theta})= -
(\lambda_1(M_N)- \rho
_{\theta})
\widehat{G}(\rho_{\theta})\widehat{G}(\lambda_1(M_N))$ and
$g_\sigma(\rho
_{\theta})={1}/{\theta},$
one gets (on $\Omega_N$)
\begin{eqnarray*}
&&\lambda_1(M_N)-\rho_{\theta}\\
&&\qquad
= \frac{(W_N)_{11}}{\sqrt{N}}+ \check{M}_{\bolds\cdot1}^*\widehat
{G}(\rho_{\theta})
\check{M}_{\bolds\cdot1}- \sigma^2 g_\sigma(\rho_{\theta}) \\
&&\qquad\quad{}
- \check{M}_{\bolds\cdot1}^*\bigl[\bigl(\lambda_1(M_N)-\rho
_{\theta
}\bigr)\widehat
{G}(\rho_{\theta})\bigl(\widehat{G}(\rho_{\theta}) -\bigl(\lambda
_1(M_N)-\rho
_{\theta
}\bigr)\widehat{G}(\rho_{\theta})\\
&&\qquad\quad\hspace*{202pt}{}\times
\widehat{G}(\lambda_1(M_N)) \bigr)\bigr]\check{M}_{\bolds\cdot1}.
\end{eqnarray*}
Finally,\vspace*{2pt} defining ${f_{\theta}(z):= \frac
{1}{\rho
_{\theta}-z}1_{\vert z \vert\leq2\sigma+ \delta}}$,
we can easily deduce from the previous equality the following identity on
$\Omega_N$:
%
\begin{eqnarray}\label{represlmax1}
&&\{ 1 + c_N + \delta_1(N)+
\delta_2(N)\} \sqrt{N} \bigl( \lambda_1(M_N)-\rho_{\theta}
\bigr) \nonumber\\[-8pt]\\[-8pt]
&&\qquad = (W_N)_{11} +
\sqrt{\frac{N}{N-1}}d_N + \sqrt{\frac{N}{N-1}} \delta_3(N),\nonumber
\end{eqnarray}
where
\begin{eqnarray*}
c_N&=&\sigma^2
\mathrm{tr}_{N-1}[f_{\theta}^2(\widehat{M}_{N-1})],
\\
d_N&=&\sqrt{N-1}\bigl\{ \check{M}_{\bolds\cdot1}^*\widehat{G}(\rho
_{\theta
})1_{\Vert\widehat{M}_{N-1}\Vert\leq2 \sigma+ \delta}\check
{M}_{\bolds\cdot1}
- \sigma^2 \mathrm{tr}_{N-1}\widehat{G}(\rho_{\theta})1_{\Vert
\widehat{M}_{N-1}\Vert\leq2 \sigma+ \delta}\bigr\},
\\
\delta_1(N)&=&-\bigl(\lambda_1(M_N)-\rho_{\theta}\bigr) \check{M}_{\bolds\cdot1}^*
[\widehat
{G}(\rho_{\theta})]^2 \widehat{G}(\lambda_1(M_N))\check
{M}_{\bolds\cdot1}1_{\Omega_N},
\\
\delta_2(N)&=&\check{M}_{\bolds\cdot1}^*\bigl[\widehat{G}(\rho_{\theta
})1_{\Vert\widehat{M}_{N-1}\Vert\leq2 \sigma+ \delta}
\bigr]^2\check{M}_{\bolds\cdot1} -
\sigma^2 \mathrm{tr}_{N-1}\bigl[\widehat{G}(\rho_{\theta})1_{\Vert
\widehat{M}_{N-1}\Vert\leq2 \sigma+ \delta}\bigr]^2,
\\
\delta_3(N)&=&\sigma^2 \sqrt{N-1} \biggl\{ \mathrm{tr}_{N-1}(
f_{\theta
}(\widehat{M}_{N-1}))- \int f_{\theta}\,d\mu_{sc}\biggr\}.
\end{eqnarray*}
First
\begin{eqnarray*}
| \delta_1(N)| &\leq&|\lambda_1(M_N)-\rho_{\theta}| \Vert\check
{M}_{\bolds\cdot1}\Vert^2 \Vert\widehat{G}(\rho_{\theta})\Vert^2 \Vert\widehat
{G}(\lambda_1(M_N))\Vert1_{\Omega_N}\\
& \leq&\frac{1}{(\rho_{\theta} - 2 \sigma- 2 \delta)(\rho
_{\theta} -
2 \sigma- \delta)^2}
\frac{1}{N} \sum_{j=2}^N |(W_N)_{j1}|^2 \times|\lambda_1(M_N)-\rho
_{\theta}|,
\end{eqnarray*}
[using Lemma \ref{lem0}(v)]. By the law of large numbers
$\frac{1}{N} \sum_{j=2}^N
|(W_N)_{j1}|^2$ converges a.s. toward $\sigma^2$ and according to Theorem
\ref{ThmASCV}, $|\lambda_1(M_N)-\rho_{\theta}|$ converges a.s. to zero.
Hence $ \delta_1(N)$ converges obviously in probability toward zero.

Now, since
$f_{\theta}$ is analytic on an open set including $[-2 \sigma,
2\sigma
]$, we deduce from Theorem
\ref{BaiYao} the convergence in probability of $\delta_3(N)$ toward
zero and of $c_N $ toward ${\sigma^2 \int f_{\theta
}^2\,d\mu_{sc}=\frac{\sigma^2}{\theta^2-\sigma^2} }$.\vadjust{\goodbreak}

According to Theorem \ref{BaiSilver98} and using Lemma \ref{lem0}(v),
\begin{eqnarray*}
\mathbb{E}(\vert\delta_2(N) \vert^2)&\leq&\frac{K}{N-1} {\mathbb
E}
\bigl( \mathrm{tr}
_N \bigl[\widehat{G}(\rho_{\theta})1_{\Vert\widehat{M}_{N-1}\Vert
\leq
2 \sigma+ \delta}\bigr]^4 \bigr) \\
& \leq& \frac{K}{N-1}
{\mathbb E}
\bigl(\Vert\widehat{G}(\rho_{\theta})\Vert^4
1_{\Vert\widehat{M}_{N-1}\Vert\leq2 \sigma+ \delta}\bigr)\\
& \leq& \frac{K}{N-1}   \frac{1}{(\rho_{\theta}-2 \sigma-
\delta)^4}.
\end{eqnarray*}
The convergence in probability of $\delta_2(N)$ toward zero readily
follows by Chebyshev inequality.

Let us check that $\widehat{G}(\rho_{\theta})1_{\Vert
\widehat{M}_{N-1}\Vert\leq2 \sigma+ \delta}$ satisfies the
conditions of Theorem \ref{Silver}.

\begin{longlist}
\item $\Vert\widehat{G}(\rho_{\theta})1_{\Vert\widehat
{M}_{N-1}\Vert\leq2 \sigma+ \delta}\Vert\leq\frac{1}{\rho
_{\theta
}-2 \sigma- \delta}$ by Lemma \ref{lem0}(v).
\item
As already noticed, $\mathrm{tr}_{N-1} f_{\theta
}^2(\widehat
{M}_{N-1})$ converges in probability toward $\int f_{\theta}^2\,d\mu
_{sc}$. Since
on the event $\{\Vert\widehat{M}_{N-1}\Vert\leq2 \sigma+ \delta\}
$, with limiting probability 1,
$\mathrm{tr}_{N-1} [\widehat{G}(\rho_{\theta})1_{\Vert\widehat
{M}_{N-1}\Vert
\leq2 \sigma+ \delta}]^2$
coincides with $\mathrm{tr}_{N-1} f_{\theta}^2(\widehat{M}_{N-1})$,
it also converges in probability toward $\int f_{\theta}^2\,d\mu
_{sc}$.
\item
It is proved in Proposition 3.1 in \cite{CD} that
for any $z \in\mathbb{C}$ such that $\Im m (z) >0$,
$\frac{1}{N-1}\sum_{i=1}^{N-1} ([\widehat{G}(z)]_{ii})^2$ converges in
probability toward
$g_{\sigma}^2(z)$. The same result holds for $\frac{1}{N-1}\sum
_{i=1}^{N-1} ([\widehat{G}(z)]_{ii})^2 1_{\Vert\widehat
{M}_{N-1}\Vert
\leq2 \sigma+ \delta}$.
For any $\epsilon> 0$ and any $\alpha>0$,
\end{longlist}\vspace*{-10pt}
\begin{eqnarray*}
&&\mathbb P\Biggl( \Biggl|
\frac{1}{N-1}\sum_{i=1}^{N-1} ([\widehat{G}(\rho_{\theta})]_{ii})^2
1_{\Vert\widehat{M}_{N-1}\Vert\leq2 \sigma+ \delta} -g_{\sigma
}^2(\rho_{\theta})\Biggr| > \epsilon\Biggr)\\
&&\qquad\leq
\mathbb P\Biggl( \Biggl\vert
\frac{1}{N-1}\sum_{i=1}^{N-1} \bigl\{([\widehat{G}(\rho_{\theta
})]_{ii})^2-\bigl([\widehat{G}(\rho_{\theta}+ i \alpha)]_{ii}\bigr)^2\bigr\}
1_{\Vert\widehat{M}_{N-1}\Vert\leq2 \sigma+ \delta} \Biggr\vert> \frac
{\epsilon}{3} \Biggr) \\
&&\qquad\quad{} + \mathbb P\Biggl( \Biggl\vert
\frac{1}{N-1}\sum_{i=1}^{N-1} \bigl([\widehat{G}(\rho_{\theta}+ i
\alpha
)]_{ii}\bigr)^2 1_{\Vert\widehat{M}_{N-1}\Vert\leq2 \sigma+ \delta}
-g_{\sigma}^2(\rho_{\theta}+ i \alpha)\Biggr\vert> \frac{\epsilon}{3}
\Biggr)\\
&&\qquad\quad{} + \mathbb P\biggl( \vert g_{\sigma}^2(\rho_{\theta})
-g_{\sigma}^2(\rho_{\theta}+ i \alpha)\vert> \frac{\epsilon}{3}
\biggr).
\end{eqnarray*}
Since
\begin{eqnarray*}
&&\bigl\{([\widehat{G}(\rho_{\theta})]_{ii})^2-\bigl([\widehat
{G}(\rho
_{\theta}+ i
\alpha)]_{ii}\bigr)^2 \bigr\} 1_{\Vert\widehat{M}_{N-1}\Vert\leq2
\sigma+ \delta} \\
&&\qquad  = [\widehat{G}(\rho_{\theta})-\widehat{G}(\rho
_{\theta
}+ i \alpha)]_{ii}
[\widehat{G}(\rho_{\theta})+\widehat{G}(\rho_{\theta}+ i \alpha)]_{ii}
1_{\Vert\widehat{M}_{N-1}\Vert\leq2 \sigma+ \delta}\\
&&\qquad  = i \alpha[\widehat{G}(\rho_{\theta}) \widehat
{G}(\rho
_{\theta}+ i \alpha)]_{ii}
[\widehat{G}(\rho_{\theta})+\widehat{G}(\rho_{\theta}+ i \alpha)]_{ii}
1_{\Vert\widehat{M}_{N-1}\Vert\leq2 \sigma+ \delta},
\end{eqnarray*}
we get by using Lemma \ref{lem0}(v)
\[
\bigl\vert([\widehat{G}(\rho_{\theta})]_{ii})^2-\bigl([\widehat{G}(\rho
_{\theta}+
i \alpha)]_{ii}\bigr)^2\bigr\vert1_{\Vert\widehat{M}_{N-1}\Vert\leq2 \sigma+
\delta} \leq\frac{2 \alpha}{(\rho_{\theta}-2 \sigma- \delta)^3}.
\]
Similarly, we get that
\[
\vert g_{\sigma}^2(\rho_{\theta})
-g_{\sigma}^2(\rho_{\theta}+ i \alpha)\vert\leq\frac{2 \alpha
}{(\rho
_{\theta}-2 \sigma)^3}.
\]
Thus, choosing $\alpha$ such that $\frac{2 \alpha}{(\rho
_{\theta}-2 \sigma- \delta)^3}< \frac{\epsilon}{3}$,
we readily deduce the convergence in probability of
\[
\frac{1}{N-1}\sum_{i=1}^{N-1} ([\widehat{G}(\rho_{\theta})]_{ii})^2
1_{\Vert\widehat{M}_{N-1}\Vert\leq2 \sigma+ \delta}
\]
toward $g_{\sigma} ^2(\rho_{\theta})$.

Since $\widehat{G}(\rho_{\theta})1_{\Vert\widehat
{M}_{N-1}\Vert\leq2 \sigma+ \delta}$ and
$\check{M}_{\bolds\cdot1}$ are independent, we can deduce from
Theorem \ref{Silver} that $d_N$ converges in distribution toward a
Gaussian law with mean zero and
variance
\[
v_{\theta}:=\sigma^4
\biggl\{ \biggl(\mathbb{E}\biggl(\biggl|\frac{(W_N)_{12}}{\sigma}\biggr|^4\biggr) -1 - t/2
\biggr)
\frac{1}{\theta^2}
+ \frac{t}{2}\frac{1}{\theta^2-\sigma^2}\biggr\},
\]
where $t=4$ in the real
setting and $t=2$ in the complex one. Note that one readily verifies
that $v_{\theta}$
satisfies (\ref{defvtheta}) in Section {\ref{results}}.

Let $0< \epsilon<1$. Since $\delta_1(N)+ \delta_2(N)$
converges in probability toward zero, the probability of the event
%
\begin{eqnarray}{\label{represtilde}}
\widetilde{\Omega}_N= \Omega_N \cap\{ | \delta_1(N)+ \delta_2(N)|
\leq\epsilon\}
\end{eqnarray}
tends to 1. Now, since $c_N \geq0$ we have the following identity
on $\widetilde{\Omega}_N$:
%
\begin{equation}{\label{represlmax2}}\quad
\quad \sqrt{N} \bigl(\lambda
_1(M_N)-\rho
_{\theta
} \bigr) = \frac{1}{u_N}
\biggl\{(W_N)_{11} + \sqrt{\frac{N}{N-1}}d_N + \sqrt{\frac{N}{N-1}}
\delta_3(N)\biggr\}
\end{equation}
with $u_N:=1+c_N + \delta_1(N)+ \delta_2(N)$ converging in
distribution toward $(1- {\sigma^2}/{\theta^2} )^{-1}$.
Moreover, since $(W_N)_{11}$ and $d_N$ are independent, $(W_N)_{11}
+\break
\sqrt{N/(N-1)}d_N + \sqrt{N/(N-1)} \delta_3(N)$ converges in
distribution toward the convolution of $\mu$ and a Gaussian
distribution $\mathcal{N}(0, v_{\theta})$.

Finally, we can conclude that $\sqrt{N} (1- {\sigma^2}/{\theta^2}
)^{-1} (\lambda_1(M_N)-\rho_{\theta})$ converges in distribution
toward $\mu\ast\mathcal{N}(0,
v_{\theta})$.

\subsection{\texorpdfstring{Proof of Theorem
\protect\ref{ThmFluctuationsUniversalite}}{Proof of Theorem 2.4}}
\label{sectiongoe}
As before, $\theta$ is assumed to be $>\sigma$. In
Theorem~\ref{ThmFluctuationsUniversalite}, we consider the real
Deformed models and claim that the full deformation $A_N$ defined by
$(A_N)_{ij} = {\theta}/{N}$ exhibits universality of the Gaussian
fluctuations of the largest eigenvalue around $\rho_{\theta}$. As
already stated, the analogue of this result holds in the complex
setting. This is one of the conclusions of the work \cite{FePe}
which also partly solves the real case (we recall to the reader that
all the results of \cite{FePe} readily extend to the framework of
Theorem \ref{ThmFluctuationsUniversalite} calling on \cite{Ru}). In
order to explain this more precisely, let us summarize the main
arguments developed by \cite{FePe} in the complex setting. First, it
is shown that the universality of the fluctuations follows from the
universality of limits of expectations of traces of suitable high
powers of any Deformed Wigner matrices (the powers are of the order
of $\sqrt
N$). Second (this is the main part of the work \cite{FePe}), to
handle such expectations, the authors perform a combinatorial method
inspired by \cite{So} and then deduce that in the large limit $N \to
\infty$, the previous expectations behave as in the Gaussian case. The
last step of the analysis calls on the investigations of \cite{Pe} on
the Deformed GUE which allow to identify the value of these limits.

Actually, the combinatorial arguments also work in the real setting
(see in particular Section 2 in \cite{FePe}) and reduce the
universality problem to the knowledge of the Deformed GOE. Thus, to
get the result of Theorem \ref{ThmFluctuationsUniversalite}, it
suffices to prove (using the orthogonal invariance of the GOE) the
following limit.

\begin{proposition}\label{expTL}
Call $L_{\theta}$ the Laplace transform of the law $\mathcal N (0, 2
\sigma_{\theta} ^2)$. Let~$M_N^G$ be the Deformed GOE
with $A_N=\operatorname{diag}(\theta, 0,\ldots,0)$ and $\theta>\sigma$.

For any $t $ in $[0,\rho_{\theta}[$,
%
\begin{equation} \label{limitetrace} \lim_N
\mathbb E\bigl[\operatorname{Tr} ({M_N^G}/{\rho_{\theta}})^{2
[t \sqrt N] }\bigr] = L_{\theta}( {2t}/{\rho_{\theta}}).
\end{equation}
\end{proposition}

The starting point of our computations is the following
result which states that the previous expectation only involves (as
$N \to\infty$) the rescaled largest eigenvalue of the Deformed GOE
\[
\xi_1 ^G = \sqrt N
\bigl( \lambda_1(M_N^G)-\rho_{\theta} \bigr).
\]

\begin{lemma} For any $t >0$,
%
\begin{equation} \label{equivtrace}
\mathbb E\bigl[\operatorname{Tr} ({M_N^G}/{\rho_{\theta}}
)^{2[t \sqrt N ]}\bigr] = \mathbb E\bigl[\exp({2t \xi_1
^G}/{\rho
_{\theta}} )   1 _
{\vert\xi_1 ^G\vert\leq N^{1/6}} \bigr] \bigl(1+o(1)\bigr).
\end{equation}
\end{lemma}

This formula does not appear explicitly in \cite{FePe} but all the
arguments needed for its justification can be found in it (actually one
can show that the formula holds for any Deformed Wigner model $M_N$
satisfying the assumptions of Theorem \ref
{ThmFluctuationsUniversalite}). We will not give the proof and refer the
reader to Section 2 in \cite{FePe}.

Hence, to derive Proposition \ref{expTL}, it remains to show the next
lemma on $\xi_1 ^G$.

\begin{lemma}\label{TL} For any $t $ in $[0,2[$,
%
\begin{equation} \label{limiteTLGOE}
\lim_N   \mathbb E\bigl[\exp
(t \xi_1^G )   1 _
{\vert\xi_1 ^G\vert\leq N^{1/6}} \bigr] = L_{\theta}( t).
\end{equation}
\end{lemma}

\begin{pf} Observe first that it is
enough to show that
%
\begin{equation} \label{limiteTLGOEbis}
\lim_N   \mathbb E \bigl[\exp
(t \xi_1^G )   1 _
{\vert\xi_1 ^G\vert\leq N^{1/6}}   1 _{\widetilde{\Omega}_N}
\bigr] = L_{\theta}( t),
\end{equation}
where the event
${\widetilde{\Omega}_N}$ was defined above by (\ref{represtilde})
choosing $ \delta>0 $ smaller than $\min\{ \frac{\rho_\theta
-2 \sigma}{4}; \frac{1}{3} \int\frac{1}{\rho_\theta-x}\, d
\mu_{sc}(x) \}$. Indeed, by the Cauchy--Schwarz inequality,
\[
\bigl( \mathbb E \bigl[ \exp(t \xi_1^G )   1 _
{\vert\xi_1 ^G\vert\leq N^{1/6}}   1 _{^c\widetilde{\Omega}_N}
\bigr] \bigr) ^2 \leq\mathbb E \bigl[ \exp(2t \xi_1^G
)   1 _ {\vert\xi_1 ^G\vert\leq N^{1/6}} \bigr] \times
\mathbb P (^c\widetilde{\Omega}_N).
\]
The previous right-hand side is negligible as $N \to\infty$ since the
probability vanishes and the expectation is bounded since \cite{FePe}
proved that the left-hand side of~(\ref{equivtrace}) is bounded, too.

The occurrence of the event ${\widetilde{\Omega}_N}$ allows to make
use of the relevant representation ({\ref{represlmax2}}) of $\xi_1
^G$ obtained in the
previous Section \ref{sectionThmfluctu}:
%
\begin{equation}\label{represlmaxGOE}
\xi_1 ^G = \frac{1}{u_N^G}
\biggl\{(W_N^G)_{11} + \sqrt{\frac{N}{N-1}}d_N^G + \sqrt{\frac{N}{N-1}}
\delta_3^G(N)\biggr\}.
\end{equation}
Second, by Fubini's theorem one can check that
%
\begin{eqnarray}{\label{reecritureTLGOE}}\qquad
&&\mathbb E \bigl[\exp(t \xi_1^G )   1 _ {\vert\xi_1
^G\vert\leq N^{1/6}}   1 _{\widetilde{\Omega}_N}
\bigr]\nonumber\\[-8pt]\\[-8pt]
&&\qquad = \int
_{\mathbb R} e^x   \mathbb P [ \{ t\xi_1^G \geq x \} \cap
\widetilde{\Omega}_N \cap\{ \vert\xi_1 ^G\vert\leq N^{1/6} \}
] \, dx.\nonumber
\end{eqnarray}
By Theorem \ref{ThmFluctuations},
%
\begin{equation}\label{fluct}\mathbb P [ \{ t\xi_1^G \geq x \}
\cap
\widetilde{\Omega}_N \cap\{ \vert\xi_1 ^G\vert\leq N^{1/6} \}
] \mathop{\longrightarrow}\limits_{N \rightarrow+ \infty}
\mathbb P
[ t \mathcal{N} \geq x ],
\end{equation}
where $\mathcal{N}$ is a
centered Gaussian variable with variance $2 \sigma_\theta^2$. We
want to deduce~(\ref{limiteTLGOEbis}) from (\ref{reecritureTLGOE})
and (\ref{fluct}) by the dominated convergence theorem. Thus, we are
going to prove that there exists a function $h$ such that for $N$
large enough and for any $x$,
\[
\mathbb P [ \{ t\xi_1^G \geq x \}
\cap\widetilde{\Omega}_N \cap\{ \vert\xi_1 ^G\vert\leq
N^{1/6} \} ] \leq h(x)
\]
with $\int_{\mathbb R} e^x   h(x)
\, dx < + \infty.$ Note that for $x \leq0$, the result is obvious
setting $h(x)=1$. Let $x$ be nonnegative. We shall improve the
general analysis made in the previous Section \ref{sectionThmfluctu}
thanks to the particular Gaussian setting considered here. For all
$N$ large enough,
\begin{eqnarray*}
& & {\mathbb P [ \{ t\xi_1^G \geq x \} \cap
\widetilde{\Omega}_N \cap\{ \vert\xi_1 ^G\vert\leq N^{1/6} \}
]}\\
&&\qquad \leq \mathbb P \biggl[
(W_N^G)_{11} \geq\frac{x (1- \varepsilon)}{3t} \biggr]
+ \mathbb P \biggl[
\sqrt{\frac{N}{N-1}}d_N^G \geq\frac{x (1- \varepsilon)}{3t}
\biggr]\\
&&\qquad\quad{}  +   \mathbb P \biggl[
\sqrt{\frac{N}{N-1}} \delta_3^G(N) \geq\frac{x (1- \varepsilon)}{3t}
\biggr]\\
&&\qquad =  J_N^{(1)}(x) +J_N^{(2)}(x)+ J_N^{(3)}(x).
\end{eqnarray*}
$J_N^{(1)}(x)=J^{(1)}(x)$ does not depend on $N$ and we
have $\int_{\mathbb R} e^x   J^{(1)}(x) \, dx = \mathbb{E}[
\exp(\frac{3t}{1-\varepsilon} (W_N^G)_{11} ) ]<
+ \infty.$ Besides, one can easily see that the choice of $\delta$
insures that $J_N^{(3)}(x)=0.$ By the Chebyshev inequality,
we have
\[
J_N^{(2)}(x) \leq\exp\bigl(-6x(1- \varepsilon)/{3t}\bigr) \mathbb{E}
(\mathcal{E}),
\]
where $\mathcal{E}= \mathcal{E}' \mathcal{E}'' $ with
\begin{eqnarray*}
\mathcal{E}' &=&\mathbb{E}\bigl[ \exp\bigl(
6 \sqrt{N}     ^t \check{M}_{\bolds\cdot1}^{G}
\widehat{G}(\rho_{\theta}) 1_{\Vert\widehat{M}^G_{N-1}\Vert\leq2
\sigma+ \delta}   \check{M}^G_{\bolds\cdot1} \bigr)\vert
 \widehat{M}^G_{N-1} \bigr],
\\
\mathcal{E}''&=&
\exp\bigl[- 6\sigma^2 \sqrt{N}
\mathrm{tr}_{N-1}\bigl(\widehat{G}(\rho_{\theta})1_{\Vert
\widehat{M}^G_{N-1}\Vert\leq2 \sigma+ \delta}\bigr)\bigr].
\end{eqnarray*}
Using the Gaussian assumptions (see \cite{Sa}, pages 90--91), one has
\begin{eqnarray*}
\mathcal{E}' & = & \det\biggl(I_{N-1} - 12\frac
{\sigma
^2}{\sqrt{N}}\widehat{G}(\rho_{\theta})1_{\Vert\widehat
{M}^G_{N-1}\Vert\leq2 \sigma+ \delta}\biggr)^{-{1}/{2}}\\
&=& \prod_{i=1}^{N-1} \bigl(1 - 12{\sigma^2}\beta_i/{\sqrt{N}}
\bigr)^{-{1}/{2}}
\end{eqnarray*}
for large enough $N$, where the $\beta_i$'s are the eigenvalues of
$\widehat{G}(\rho_{\theta})1_{\Vert\widehat{M}^G_{N-1}\Vert\leq2
\sigma+ \delta}$. Note that $0\leq\beta_i < \frac{1}{3 \delta}$
so that the last identities make sense, for instance, for $N >
\frac{16 \sigma^4}{ \delta^2}$. Hence,
\begin{eqnarray*}
\ln\mathcal{E}' \mathcal{E}''
\leq\frac{1}{2}\sum_{i=1}^{N-1} \bigl\{ -\ln\bigl(1 -
12{\sigma^2}\beta_i/{\sqrt{N}}\bigr) -
12{\sigma^2}\beta_i/{\sqrt{N}}\bigr\}.
\end{eqnarray*}
Let $\alpha> \frac{1}{2}$; using that for any $y$ in $[0,
1-\frac{1}{2\alpha}]$, we have $-\ln(1-y) -y \leq\alpha y^2$. So,
as $\beta_i < \frac{1}{3 \delta}$, we get that for $N > {16
\sigma^4}{\delta^{-2}(1-\frac{1}{2\alpha})^{-2}}$,
\[
\ln\mathcal{E}'
\mathcal{E}''\leq\frac{\alpha12^2 \sigma^4}{18 \delta^2}.
\]
Thus,
there is some constant $C_{\alpha, \sigma, \delta}$ such that
$\mathcal{E}' \mathcal{E}''\leq C_{\alpha, \sigma, \delta}$ and
\[
J^{(2)}_N(x)\leq C_{\alpha, \sigma, \delta} \exp\bigl(-2x {(1-
\varepsilon)}/{t}\bigr).
\]
Now, for $0<t <2 (1- \varepsilon)$, $\int_0^{+ \infty} \exp
(x-2x {(1- \varepsilon)}/{t})\,dx < \infty$. The
proof is complete.
\end{pf}

\begin{appendix}\label{app}
\section*{Appendix: By J. Baik and J. Silverstein}

This \hyperref[app]{Appendix} presents the proof by J. Baik and J. Silverstein of the
CLT (given by Theorem \ref{Silver})
needed in the previous section for the proof of Theorem
\ref{ThmFluctuations}. Their proof is based on a
writing of the expression
\setcounter{equation}{0}
\begin{equation}{\label{variable}}
\bigl(1/\sqrt N\bigr) (Y_N^*BY_N-\operatorname{Tr} B )
\end{equation}
as a sum of
martingale differences, and uses the following CLT.

\setcounter{theorem}{0}
\begin{theorem}[(Theorem 35.12 of \protect\cite{Bil})] \label{Theo-Bil} For each
$N$, let $Z_{N1}, \ldots, Z_{Nr_N}$ be a real martingale difference sequence
with respect to the increasing $\sigma$-field $\{ \mathcal F _{N,j} \}$
having second moments. If, as $N \to\infty$,
%
\begin{equation}{\label{Condition1}}
\sum_{j=1}^{r_N} \mathbb E ( Z_{Nj}^2 | \mathcal F _{N,j-1})
\stackrel{P}{\longrightarrow} v^2,
\end{equation}
where $v^2$ is a positive constant, and for each $\epsilon>0$,
%
\begin{equation}{\label{Condition2}}
\sum_{j=1}^{r_N} \mathbb E \bigl( Z_{Nj}^2   1 _{| Z_{Nj}| \geq\epsilon}
\bigr)
\rightarrow  0,
\end{equation}
then
\[
\sum_{j=1}^{r_N} Z_{Nj} \stackrel{\mathcal L }{\longrightarrow}
\mathcal N (0,v^2).
\]
\end{theorem}

\begin{pf*}{Proof of Theorem \protect\ref{Silver}} First, one can write
(\ref{variable}) as a sum of martingale differences:
\begin{eqnarray*}
&&\bigl(1/\sqrt N\bigr)(Y_N^*BY_N-\operatorname{Tr} B)\\
&&\qquad=\bigl(1/\sqrt N\bigr) \sum_{i=1}^N
\Biggl((|y_i|^2-1)b_{ii}
+\bar y_i\sum_{j<i}y_jb_{ij}+\bar y_i\sum_{j>i}y_jb_{ij}\Biggr)
\\
&&\qquad=\bigl(1/\sqrt N\bigr) \sum_{i=1}^N
\Biggl((|y_i|^2-1)b_{ii}
+\bar y_i\sum_{j<i}y_jb_{ij}+y_i\sum_{j<i}\bar y_jb_{ji}\Biggr)
=\sum_{i=1}^NZ_i,
\end{eqnarray*}
where
\[
Z_i=Z_{Ni}=\bigl(1/\sqrt N\bigr)
\Biggl((|y_i|^2-1)b_{ii}
+\bar y_i\sum_{j<i}y_jb_{ij}+y_i\sum_{j<i}\bar y_j\bar b_{ij}\Biggr).
\]
Let $\mathcal F _{N,i}$ (resp., $\mathcal F _{N,0}$) be the
$\sigma$-field generated by $y_1,\ldots,y_i$ and $B$ (resp., by~$B$).
Let also $\mathbb{E}_i(\cdot)$ denote conditional expectation with
respect to $\mathcal F _{N,i}$. It is clear that~$Z_i$ is measurable
with respect to $\mathcal F _{N,i}$ and satisfies $\mathbb E_{i-1} (
Z_{i})=0$.

We will show the conditions of Theorem \ref{Theo-Bil} are met.

To verify the Lindeberg condition (\ref{Condition2}), we need to show
this property
is closed under addition. This will follow from the following fact.
For random variables $X_1$, $X_2$, and positive $\epsilon$,
%
\begin{equation}\label{addition}\hspace*{30pt}
\mathbb{E}\bigl(|X_1+X_2|^2   1_{|X_1+X_2|\ge\epsilon}\bigr) \leq4 \bigl(
\mathbb{E} \bigl( |X_1|^2   1_{|X_1|\ge\epsilon/2}\bigr)+\mathbb{E}\bigl(|X_2|^2
  1_{|X_2|\ge\epsilon/2} \bigr) \bigr).
\end{equation}
Indeed, we have
\begin{eqnarray*}
\mathbb{E}\bigl(|X_1|^2   1_{|X_1+X_2|\ge\epsilon}\bigr) & \leq&
\mathbb{E}\bigl(|X_1|^2   1_{(|X_1|\ge\epsilon/2)}\bigr) +\mathbb{E}\bigl(|X_1|^2
  1_{(|X_1|<\epsilon/2,|X_2|\ge\epsilon/2)}\bigr)\\
& \leq&
\mathbb{E}\bigl(|X_1|^2
1_{(|X_1|\ge\epsilon/2)}\bigr)+(\epsilon^2/4)\mathbb
P(|X_2|\ge\epsilon/2)\\
& \leq& \mathbb{E}\bigl(|X_1|^2
1_{(|X_1|\ge\epsilon/2)}\bigr)+\mathbb{E}\bigl(|X_2|^2
1_{(|X_2|\ge\epsilon/2)}\bigr).
\end{eqnarray*}
The same bound starting with $X_2$ leads to (\ref{addition}).

Write $Z_i=X_1^i+X_2^i$, with $X_1^i=(1/\sqrt
N)(|y_i|^2-1)b_{ii}$. Then for $\epsilon>0$,
%
\begin{equation}\label{terme1}
\sum_{i=1}^N\mathbb{E}\bigl(|X^i_1|^2   1_{(|X_1^i|\ge\epsilon)}\bigr)\leq
a^2\mathbb{E}\bigl((|y_1|^2-1)^2   1_{(||x_1|^2-1|\ge\sqrt N\epsilon
/a)}\bigr)\to0
\end{equation}
as $N \to\infty$, by the dominated convergence theorem.

We have
\begin{eqnarray*}
\mathbb{E}\Biggl|\sum_{j<i}y_jb_{ij}\Biggr|^4 & = & \mathbb{E}\Biggl(|y_1|^4\sum
_{j<i}|b_{ij}|^4\Biggr)
+2 \mathbb{E} \Biggl(\sum_* |b_{ij_1}|^2|b_{ij_2}|^2\Biggr)\\
&&{} + \mathbb
{E}\Biggl(|y_1^2|^2\sum_*
b_{ij_1}^2\bar b_{ij_2}^2 \Biggr)\\
& \leq&
\mathbb{E}|y_1|^4 \mathbb{E}\biggl[\max_j(B^2)_{jj}(B^2)_{ii}\biggr]+(2+\mathbb
{E}|y_1^2|^2)
\mathbb{E}[(B^2)_{ii}^2] \\
& \leq& a^4 [ \mathbb{E}|y_1|^4 + 2+\mathbb{E}|y_1^2|^2 ],
\end{eqnarray*}
where the sum $\sum\limits_*$ is over $\{j_1<i,   j_2<i,   j_1
\not= j_2 \}$.
Therefore $\mathbb{E}|X_2^i|^4=o(N^{-1})$ so that for any $\epsilon>0$,
%
\begin{equation}\label{terme2}
\sum_{i=1}^N\mathbb{E}\bigl(|X_2^i|^2   1_{(|X_2^i|\ge\epsilon)}\bigr)\leq
(1/\epsilon^2)
\sum_{i=1}^N\mathbb{E}|X_2^i|^4
\to0 \qquad\mbox{as } N\to\infty.
\end{equation}
Thus, by (\ref{terme1}), (\ref{terme2}) and (\ref{addition}), $\{
Z_i\}$
satisfies (\ref{Condition2}).

Now, we shall verify condition (\ref{Condition1}). We have
%
\begin{eqnarray}\label{sum-cond1}\qquad\quad
&&\sum_{i=1}^N\mathbb{E}_{i-1}Z_i^2\nonumber\\[-2pt]
&&\qquad
= (1/N) \sum_{i=1}^N \Biggl\{ (\mathbb{E}|y_1|^4-1)b_{ii}^2
+\mathbb{E}\bar
y_1^2\Biggl(\sum_{j<i}y_jb_{ij}\Biggr)^2\nonumber\\[-9pt]\\[-9pt]
&&\qquad\quad\hspace*{49.3pt} {}+ \mathbb{E}
y_1^2\Biggl(\sum_{j<i}\bar y_j\bar b_{ij}\Biggr)^2
+2\mathbb{E}(|y_1|^2\bar
y_1)b_{ii}\sum_{j<i}y_jb_{ij}\nonumber\\[-2pt]
&&\qquad\quad\hspace*{49.3pt}{} +
2\mathbb{E}(|y_1|^2y_1)b_{ii}\sum_{j<i}\bar y_j\bar b_{ij}
+
2\Biggl(\sum_{j<i}y_jb_{ij}\Biggr)\sum_{j<i}\bar y_j\bar b_{ij} \Biggr\}.
\nonumber
\end{eqnarray}
Let $B_L$ denote the strictly lower triangular part of $B$.
We have
\[
\mathbb{E}\Biggl[(1/N)\sum_{i=1}^N b_{ii}\sum_{j<i}y_jb_{ij}\Biggr]=0
\]
and using Cauchy--Schwarz,
\begin{eqnarray*}
\mathbb{E}\Biggl|(1/N)\sum_{i=1}^N b_{ii}\sum_{j<i}y_jb_{ij}\Biggr|^2&
=&
\mathbb{E}\Biggl|(1/N)\sum_{j=1}^{N-1}
y_j\sum_{i>j}b_{ii}b_{ij}\Biggr|^2\\
&=&  (1/N^2) \mathbb{E}
\Biggl(\sum_{j=1}^{N-1}\sum_{i>j}b_{ii}b_{ij}\sum_{\underline i>j}
b_{\underline i \underline i}\overline b_{\underline ij}\Biggr)\\
&=&(1/N^2) \mathbb{E} \Biggl( \sum_{i\underline i}b_{ii}b_{\underline
i \underline i} (B_LB_L^*)_{i\underline i}\Biggr)\\
&\leq&
\mathbb{E} \Biggl[ \biggl(\max_ib_{ii}\biggr)^2(1/N) \Biggl(\sum_{i\underline
i}|(B_LB_L^*)_{i\underline i}|^2\Biggr)^{1/2} \Biggr]\\
&=& \mathbb{E} \Biggl[ \biggl(\max_ib_{ii}\biggr)^2(1/N)
\operatorname{Tr}((B_LB_L^*)^2)^{1/2} \Biggr]\\
&\leq& \mathbb{E}
\biggl[ \biggl(\max_ib_{ii}\biggr)^2\bigl(1/\sqrt N\bigr)\|B_L\|^2 \biggr].
\end{eqnarray*}
We apply the following bound (due to Mathias; see \cite{Mt}):
$\|B_L\|\leq\gamma_N\|B\|$ where $\gamma_N=O(\ln N)$, and the bound
$\Vert B\Vert  \leq a$ to conclude that
\[
(1/N) \sum_{i=1}^N b_{ii}\sum_{j<i}y_jb_{ij} \stackrel
{P}{\longrightarrow
} 0.
\]
Then (recall that $\mathbb{E} y_1^2=0$ when $y_1$ is complex), (\ref
{sum-cond1}) can be written as
%
\begin{eqnarray}{\label{sumbis-cond1}}\qquad\quad
\sum_{i=1}^N\mathbb{E}_{i-1}Z_i^2 & = & (1/N)\sum_{i=1}^N
\Biggl[(\mathbb{E}|y_1|^4-1)b_{ii}^2\nonumber\\
&&\hspace*{48.8pt}{}
+t\Biggl(\sum_{j<i}y_jb_{ij}\Biggr)\Biggl(\sum_{j<i}\bar y_j\bar
b_{ij}\Biggr)\Biggr]+o_{P}(1)\\
& = & (1/N)\sum_{i=1}^N
(\mathbb{E}|y_1|^4-1)b_{ii}^2+t(1/N)Y_N^*B_L^*B_LY_N+o_{P}(1),\nonumber
\end{eqnarray}
where $t=4$ when $y_1$ is real, and is 2 when $y_1$ is complex.

Besides, from Lemma 2.7 in \cite{BS1} (recalled in Theorem
\ref{BaiSilver98}) we have
\begin{eqnarray*}
\mathbb{E} \bigl|(1/N)\bigl(Y_N^*B_L^*B_L Y_N-\operatorname{Tr} (B_L^*B_L) \bigr)
\bigr|^2 & \leq&
(1/N^2) \mathbb{E} (\operatorname{Tr}(B_L^*B_L)^2 ) \\
& \leq& K \mathbb{E}\|B\|^4 \frac{\ln^4N} N\to0
\end{eqnarray*}
as $N\to\infty$. So, as
\[
\operatorname{Tr} B_L^*B_L=\sum_{j<i}|b_{ij}|^2=(1/2)\Biggl(\operatorname{Tr}
B^2-\sum_ib_{ii}^2\Biggr),
\]
(\ref{sumbis-cond1}) implies that condition (\ref{Condition1}) holds
with
\[
v^2=(\mathbb{E}|y_1|^4-1-t/2)a_1^2+(t/2)a_2.
\]
Thus, by
Theorem \ref{Theo-Bil}, we deduce that $({1}/{\sqrt N})
(Y_N^*BY_N- \operatorname{Tr} B )$ converges in distribution to a
Gaussian variable with mean zero and variance $v^2$.
\end{pf*}
\end{appendix}

\section*{Acknowledgments}
The authors are very grateful to
Jack Silverstein and Jinho Baik for providing them their proof of
Theorem \ref{Silver} (which is a fundamental argument in the proof of
Theorem \ref{ThmFluctuations}) presented in the \hyperref[app]{Appendix} of
the present article.
The authors also wish to thank an anonymous referee for useful
comments which led to an improvement of this paper.


%
\printaddresses

\end{document}